\documentclass[12pt]{amsart}
\usepackage{amssymb,amscd}
\usepackage[russian]{babel}
\usepackage{graphicx}
\usepackage{pdfpages}
\usepackage{placeins}
\newtheorem{theorem}{Теорема}[section]

\newtheorem{lemma}{Лемма}[section]

\def\Ad{\operatorname{Ad}}
\def\DOLL{\operatorname{DOLL}}

\theoremstyle{definition}
\newtheorem{definition}{Определение}[section]

\theoremstyle{remark}

\ifx\pdfpagewidth\undefined \else
\fi
 \evensidemargin \oddsidemargin

\advance\hoffset by -1in
\advance\voffset by -1in

\usepackage[cp1251]{inputenc}
\author{И.~А.~Иванов-Погодаев, А.~Я.~Канель-Белов}

\thanks{Moscow Institute of Physics and Technology, Bar-Ilan University 
and College of Mathematics and Statistics, Shenzhen University, Shenzhen, 518061, China}

\address{Moscow Institute of Physics and Technology, Russia, Bar-Ilan University, Israel}
\email{ivanov.pogodaev@mail.ru}
\address{College of Mathematics and Statistics, Shenzhen University, Shenzhen, 518061, China}
\email{kanel@mccme.ru}

\title{Комплексы со свойством равномерной эллиптичности}

\begin{document}

\let \mathbf=\texttt
\tabcolsep 2pt

\begin{abstract}

Работа посвящена конструкции конечно определенной бесконечной нильполугруппы, удовлетворяющей тождеству $x^9=0$. Эта конструкция отвечает на проблему Л.~Н.~Шеврина и М.~В.~Сапира. 

В силу большого объема доказательство разделено на геометрическую, комбинаторную и финализационную части.

В первой части строится равномерно-эллиптическое пространство. Пространство называется {\it равномерно-эллиптическим}, если любые две точки $A$ и $B$ на расстоянии $D$ соединяются системой геодезических,  образующих диск ширины $\lambda\cdot D$ для некоторой глобальной константы $\lambda>0$.

Во второй части рассматриваются комбинаторные свойства комплекса. Вершины и ребра кодируются конечным числом букв, что позволяет рассмотреть полугруппу путей. Определяющие соотношения соответствуют парам эквивалентных коротких путей на комплексе. Кратчайшим путям в смысле естественной метрики будут соответствовать ненулевые слова в полугруппе. Слова, не соответвующие путям на комплексе или соответствующие некратчайшим путям, приводятся к нулю.

В третьей части проводится финализация, в частности показывается, что слово, содержащее девятую степень может быть приведено к нулю с помощью определяющих соотношений.

Настоящая работа посвящена первой части доказательства.

Данная работа была проведена с помощью Российского Научного Фонда Грант N 17-11-01377.  Первый автор является победителем конкурса ``Молодая математика России''.

%This paper is devoted to construction of finitely presented infinite nil semigroup with identity $x^9=0$. This construction answers to the problem of Lev Shevrin and Mark Sapir.

%The paper is quite long so the proof is separated into geometric, combinatorial and finalization parts. 

%In the first part we construct uniformly elliptic space. Space is called {\it uniformly elliptic} if any two points $A$ and $B$ at the distance of $D$ can be connected by the system of geodesics which form a disc with width $\lambda\cdot D$ for some global constant  $\lambda>0$.

%In the second part we study combinatorial properties of the constructed complex. Vertices and edges of this complex coded by finite number of letters so we can consider semigroup of paths. Defining relations correspond to pairs of equivalent short paths on the complex. Shortest path in sense of natural metric correspond nonzero words in the semigroup. Words which are not presented as paths on complex and words correspond to non shortest paths can be reduced to zero. 

%In the third part we make a finalization. In particular, we show that word containing ninth degree word can be reduced to zero by defining relations.

%The present paper contains first part of the proof.

%This work was carried out with the help of the Russian Science Foundation Grant N 17-11-01377. The first author is the winner of the contest `` Young Mathematics of Russia ''.

\end{abstract}

\maketitle

\medskip
\tableofcontents

\medskip

\section{Введение} \label{nachalo}

Работа посвящена построению конечно определенных нильполугрупп. Доказана следующая

\medskip
{\bf Теорема.} {\it Существует конечно определенная бесконечная нильполугруппа, удовлетворяющая тождеству $x^9=0$.}
\medskip

%Теории полугрупп посвящено ряд обзоров (см. например, \cite{Shevrin,Sapirbook}).
 Подобные построения затрагивают, с одной стороны, проблемы бернсайдовского типа, а с другой -- проблемы построения конечно определенных объектов. Поэтому следует, хотя бы кратко, упомянуть историю вопроса в этой связи.

\subsection{Проблемы бернсайдовского типа} \label{vvedenie}

Проблемы бернсайдовского типа оказали огромное влияние на развитие современной алгебры. Эта проблематика охватила большой круг вопросов,
как в теории групп, так и в смежных областях, стимулировала алгебраические исследования. Проблемам бернсайдовского типа
 посвящена обзорная статья Е.~И.~Зельманова \cite{Zelmanov}. 

 Ясно, что группы, удовлетворяющие тождеству $x^2=1$ коммутативны. Для групп с тождеством $x^3=1$ доказать локальную конечность несколько труднее (и она была доказана самим Бернсайдом \cite{Burnside}). Вопрос о локальной конечности групп с тождеством $x^4=1$ стоял открытым чуть меньше 40 лет \cite{Sanov}, а с тождеством $x^6=1$~-- свыше 50 \cite{HallM}.

Первый контрпример к неограниченной проблеме был найден Е.~С.~Голодом в 1964 году на  основе универсальной конструкции Е.~С.~Голода -- И.~Р.~Шафаревича \cite{Golod,GolodShafarevich}. Эта конструкция позволила также построить не локально конечные периодические группы неограниченной экспоненты.

В ограниченном случае, получению контрпримеров предшествовало построение бесквадратных слов над трехбуквенным алфавитом, а также бескубных слов над алфавитом из двух букв, т.е. построению бесконечных нильполугрупп с тождествами $x^2=0$ и $x^3=0$.  Вопрос о локальной  конечности групп с тождеством $x^n = 1$ был решен отрицательно в знаменитых работах П.~С.~Новикова и С.~И.~Адяна  \cite{Novikov-Adyan}, \cite{Adyan2}: для любого нечетного $n > 4381$ было доказано существование бесконечной группы с $m>1$ образующими, удовлетворяющей тождеству $x^n = 1$. Результаты Новикова-Адяна А.~И.~Мальцев рассматривал как основное событие алгебры XX века.
 Эта оценка  была улучшена до $n > 665$ С.~И.~Адяном \cite{Adyan}. Совсем недавно С.~И.~Адян в кратком анонсе сообщил, что удалось улучшить оценку до $n \ge 101$ (см.  \cite{Adyan1}).

 Работы П.~С.~Новикова и С.~И.~Адяна оказали огромное влияние на творчество И.~А.~Рипса, который в дальнейшем разработал метод канонической формы и также построил примеры бесконечных периодических групп, обладающих дополнительными свойствами.  Позднее А.~Ю.~Ольшанский предложил геометрически
наглядный вариант доказательства для нечетных $n>10^{10}$ \cite{Olshansky,Olshansky2}. Для достаточно больших четных $n$ примеры бесконечных 2-порожденных групп периода $n$ были построены независимо С.~Ивановым и И.~Лысенком. В условиях делимости на достаточно большую степень двойки, было исследовано строение свободных бернсайдовых групп \cite{Lysenok,Lysenok2, Ivanov, IvanovS}. Подробная библиография и история вопроса -- см. \cite{Adyan}.

Бернсайдова проблематика естественным образом переносится из теории групп в теорию колец. В $PI$-случае вопросы локальной конечности алгебраических алгебр решаются положительно. В ассоциативном случае
соответствующий результат был получен И.~Капланским и Д.~Левицким. Чисто комбинаторное доказательство для ассоциативного случая получается из теоремы Ширшова о высоте \cite{Shirshov1,Shirshov2}.
 Теореме Ширшова о высоте посвящена обширная библиография (см., например, работы \cite{BelovRowenShirshov} ,\cite{Kem09},\cite{BelovHeightObzor}, \cite{BelovKharitonov}).
  Для $PI$-алгебр Ли соответствующий результат для нулевой характеристики был получен  А.~И.~Кострикиным. Кроме того, в характеристике $p$ им вместе с Е.~И.~Зельмановым показана локальная нильпотентность алгебр Ли с тождеством $y\circ \Ad(x)^n=0$. В общем случае результат был получен Е.~И.~Зельмановым. Это позволило решить так называемую ``ослабленную проблему Бернсайда'', то есть доказать наличие максимальной конечной группы среди $k$-порожденных групп, удовлетворяющих тождеству $x^n=1$.   Подробная библиография по этому вопросу изложена в монографии  \cite{Kostrikin}.

Проблемы бернсайдовского типа для полугрупп рассматривались в  монографии М.~В.~Сапира, с участием М.~В.~Волкова и В.~С.~Губы \cite{Sapirbook}.

\subsection{Проблемы конечной определенности}

Все имеющиеся примеры бесконечных периодических групп бесконечно определены. Чрезвычайно глубоким и вдохновляющим является следующий открытый вопрос (входящий список основных алгебраических проблем в теории групп):

\medskip

{\bf Вопрос.}\
{\it Существует ли конечно определенная бесконечная периодическая группа?
}
\medskip

Известна классическая теорема Хигмана о вложении рекурсивно определенных групп в конечно определенные.

В работе А.~Ю.~Ольшанского и М.~В.~Сапира \cite{OlshanskiiSapir} была построена конечно определенная группа являющаяся расширением конечно порожденной бесконечной периодической группы с помощью циклической.

На проблематику, связанную с построением разного рода
 экзотических объектов с помощью конечного числа определяющих соотношений обратил
внимание В.~Н.~Латышев ~\cite{Latyshev}. Им же была поставлена проблема существования конечно
определенного нилькольца~\cite{Dnestrovsk}.

\medskip
{\bf Вопрос (В.~Н.~Латышев).} {\it Существует ли конечно определенное бесконечномерное нилькольцо?
}
\medskip

В качестве продвижения в решении этого вопроса можно
 рассматривать результаты Г.~П.~Кукина, В.~Я.~Беляева о вложениях
 рекурсивно определенных объектов в конечно определенные \cite{Kukin,Belyaev}. В.~А.~Уфнаровским был построен пример конечно определенной алгебры промежуточного роста \cite{Ufnar1}. В работе В.~В.~Щиголева \cite{Schigolev} была изучена связь между понятиями ниль и нильпотентности конечно определённых алгебр в зависимости от количества определяющих соотношений и порождающих. Также построен пример алгебр с малым количеством определяющих соотношений, у которых все слова длины два нильпотентны. Кроме того, построению конечно определенных объектов и автоматному подходу в алгебраических структурах посвящен ряд других результатов \cite{IuduDisser,iva-fpm, PiontkZDiv, Piont, Uf2,Ufnar1, Belov,BK,Sapir}. 

\medskip

Фундаментальную проблему существования конечно определенной нильполугруппы поставили Л.~Н.~Шеврин и М.~В.~Сапир в Свердловской Тетради (3.61б) \cite{Sverdlovsk},  а также вопрос 3.8 в \cite{Obzor}.
%Философия проблемы -- в возможности добиться локальными средствами глобального эффекта (в данном случае, свойства нильности).

\medskip
{\bf Вопрос (Л.~Н.~Шеврин, М.~В.~Сапир).} {\it Существует ли конечно определенная бесконечная нильполугруппа?
}
\medskip

Авторам представляется, что полугрупповой вопрос предшествует кольцевому, а кольцевой -- групповому. Этот вопрос привлекал внимание авторов в течение многих лет. Ряд результатов, таких как конструкция конечно определенной полугруппы с нецелой размерностью Гель\-фан\-да--Ки\-рил\-ло\-ва, построение алгебр с конечным базисом Гребнера но неразрешимой проблемой делителей нуля возникли из работы над этой проблемой. \cite{BI,BI2,ivamalev}. Также был построен пример конечно определенной полугруппы, содержащий ненильпотентый ниль-идеал \cite{ivamalevsapir}.

\begin{theorem}[\cite{ivadis}]
Существует конечно определенная полугруппа H, с множеством слов $G$ от образующих, удовлетворяющая
следующим свойствам:

\begin{itemize}
  \item Существует идеал, содержащий бесконечное множество элементов $I = LH$,
где $L$~-- буква в $H$;
  \item Если слово $A\in G$ представляется в виде $A
= XYYZ$, где $X,Y,Z\in G$, тогда $LA = 0$.
\end{itemize}

\end{theorem}

{\bf Благодарности.} Авторы признательны руководителям семинара <<Теория колец>> на кафедре Высшей Алгебры механико-математического факультета МГУ В.~Н.~Латышеву и А.~В.~Михалеву за полезные обсуждения и постоянное, в течение ряда лет, внимание к работе. Мы также благодарны  И.~А.~Рипсу, Л.~Н.~Шеврину, А.~Х.~Шеню, Н.~К.~Верещагину, А.~Эршлер за полезные обсуждения, Ф.~Дюранду, Ц.~Селле, Л.~А.~Бокутю, Ю. ~Чэну, Т.~Фернику за поддержку в участии на конференциях,  А.~С.~Малистову за помощь в оформлении статьи. Особую благодарность мы выражаем А.~Л. ~Семенову за полезные советы и внимание к работе.

Результат о построении конечно определенной нильполугруппы докладывался на следующих конференциях и семинарах: \cite{China,Batumi,mipt,Marseille,Poland,tula,ippi,alglog,Paris,Sofia,Mexico}.

\section{Конечная определенность, связь с информатикой, геометрические методы}      \label{ScFnDefInfGen}

\subsection{Построение конечно определенных объектов и системы конечных автоматов}

При построении разного рода объектов трудности вызывает контроль над следствиями из вводимых соотношений, в особенности тогда, когда
следует доказать, что некое соотношение не является следствием заданных.
Зачастую используются три метода контроля над соотношениями:

\begin{enumerate}

\item Базис Гребнера и бриллиантовая лемма;
\item Теория малых сокращений;
\item Реализация машины Тьюринга или машины Минского.

\end{enumerate}

\medskip

В конечно определенном случае  вопросы, связанные с построением объектов, обладающих заданными свойствами, сильно усложняются и наибольшее значение приобретает третий метод.
При этом буква интерпретируется как состояние конечного автомата, а слово -- как цепочка взаимодействующих конечных автоматов. Если число соотношений конечно, то это взаимодействие {\it локально} и мы получаем связь с задачами самоорганизации, типа задачи Майхилла о стрелк\'ах и т.д. На этом пути была решена задача
о построении конечно определенных полугрупп с рекурсивной размерностью Гельфанда-Кириллова \cite{ivadis}. Отметим, что подобная техника довольно громоздка для построений, требующих малого роста. Например, не удалось пока построить конечно определенную полугруппу с размерностью Гельфанда-Кириллова равной $2.5$.

 На этом же пути был получен ответ на известный открытый вопрос, поставленный В.~Н.~Латышевым -- была построена  алгебра с неразрешимой проблемой делителей нуля и конечным базисом Гр\"еб\-не\-ра \cite{ivadis}, а также алгебра с неразрешимой проблемой нильпотентности и конечным базисом Гр\"еб\-не\-ра \cite{ivamalev}. Отметим, что для автоматных мономиальных алгебр (в частности, конечно определенных) а также {\it алгебр с ограниченной переработкой} аналогичный вопрос (также поставленный В.~Н.~Латышевым) решается положительно \cite{BBL,Belov,Iudu}.

\medskip

Задача о построении конечно определенной бесконечной нильполугруппы имеет интерпретацию в этих терминах. Рассмотрим цепочку локально взаимодействующих конечных автоматов. У них есть цвета корпусов. Если автомат объявляет себя нулем (совершает самоубийство), то вся цепочка погибает. Можно ли добиться того, чтобы преобразования были обратимы (если $u=v$, то $v=u$),  при этом существовали сколь угодно длинные живые цепочки, и чтобы любая цепочка, у которой несколько раз подряд повторились цвета корпусов, погибала. Мы задаем локальный закон взаимодействия, а наш враг -- задает внутренние состояния автоматов. Требуется достичь адекватного поведения цепочки автоматов.

Хотя решение проблемы построения бесконечной нильполугруппы было достигнуто геометрическими методами, данная интерпретация демонстрирует связь с самоорганизующимися системами и может быть интересной с точки зрения получения результата в этой области.

\subsection{Схема доказательства. Мозаики}

Пусть $W$~-- бесквадратное слово над алфавитом из трех букв. Если каждое его неподслово (т.е. антислово) объявить нулем, то естественно возникающая полугруппа слов обладает тождеством $x^2=0$. Однако естественная конструкция, связанная с заданием множества нулевых слов как подслов слов из некоторого семейства в конечно определенном случае работает плохо.

Невозможно задать конечно определенную ненильпотентную нильполугруппу только {\it мономиальными} соотношениями, т.е. соотношениями типа $v=0$. Ибо если есть конечный список запретов и бесконечное слово без запрещенных подслов, то есть и бесконечное {\em периодическое} слово также без запрещенных подслов.

Действительно, пусть $W$ -- бесконечное слово без запрещенных подслов. Разобьем его на подслова длины $l$, где $l$ превосходит длину каждого из запрещенных подслов. В силу бесконечности $W$ и конечности алфавита, в $W$ есть подслово вида $ABA$, где длина $A$ равна $l$. Тогда бесконечное периодическое слово с периодом $AB$ не содержит запрещенных подслов.

\medskip

Что если от одномерных расположений (слов) перейти к двумерным?
Известно, что  существуют конечные наборы многоугольников (плиток) которыми можно
замостить плоскость лишь непериодически. Впервые такой набор был построен Робертом Бергером \cite{Berger}. В
дальнейшем были построены более простые примеры, например, Рафаэлем Робинсоном \cite{Robinson}. Широко известна также мозаика Пенроуза.
Итак, имеются контактные правила, для которых: \begin{itemize}
  \item Существует замощение всей плоскости, удовлетворяющее запретам.
  \item Однако таких периодических замощений не существует.
\end{itemize}
Поэтому если бы можно было бы умножать слева-справа-сверху-снизу, то такого рода объекты можно было бы построить.

Но как придать всему этому смысл? Будем интерпретировать элементы полугруппы как пути на мозаике (дальнейший анализ показывает, что удобно иметь дело с кратчайшими путями -- иначе можно много раз проходить один и тот же цикл). Буквы кодируют плитки и переходы между ними. Если локальный непорядок (два символа плитки без символа перехода между ними и т.д.) -- то произведение ноль. Кроме того, если локальный участок не располагается на мозаике, или не располагается как участок кратчайшего пути -- то он также нулевой.

Если же любую пару узлов, соединяемую кратчайшим путем с кодом $s_1$ можно соединить
кратчайшим путем с кодом $s_2$ и наоборот, то $s_1=s_2$.

\medskip

Итак, пусть есть периодическое слово $U=W^n$. Начинаем добавлять клетки к слову $U$, локально перебрасывая пути. Оно окружается мозаикой. Поскольку $U$ периодично, то не может быть расположено на нашей мозаике. Поэтому в какой-то момент вставлять клетки не получится и мы доберемся до локального участка, несовместимого с мозаикой. Тем самым устанавливается равенство слова $U=W^n$ нулю.

Таким образом, возникает мозаика со своей глобальной структурой, которая и обеспечивает апериодичность. Локальные правила задают эту самую структуру, и путь, ``перекидывание'' которого задает область на мозаике. Возникают три языка: геометрический язык, описывающий глобальное поведение комплекса, комбинаторно геометрический язык контактов (локальных правил) и полугрупповой язык соотношений -- переброски путей. Для решения задачи следует научиться переводить с одного языка на другой и, главное, обеспечить саму эту возможность.

\medskip

Возможность перевода с языка контактов на глобальный язык связана с утверждениями типа теоремы Х.~Гудмана-Штраусса \cite{GoodmanStrauss}.  Теорема Х.~Гудмана-Штраусса утверждает, что любую плоскую мозаику, связанную с подстановочной системой (плоской $\DOLL$-системой) можно задать локальными правилами. Наш подход, использующий реберную структуру, ближе к подходу Тома Ферника и Николя Бедериде \cite{Fernique}.

\subsection{Апериодические мозаики. ``Демо-версия'' доказательства}

\subsubsection{Апериодические замощения}
Наличие апериодических мозаик связано со следующим обстоятельством. Рассмотрим работу машины Тьюринга на клетчатой ленте.  Движение головки и обмен сигналами можно перевести на язык укладки плиток на фазовой диаграмме ``пространство-время''. Из этого выводится алгоритмическая неразрешимость проблемы дополнения заданного набора плиток до замощения плоскости.

Алгоритмическая неразрешимость проблемы замощения, когда изначально ничего не положено, доказывается технически сложнее. Здесь ``пространство'' и ``время'' перемешиваются. Строятся замощения, в некоторых местах которых вынуждается один шаг машины Тьюринга, в других -- два шага и т.д. Линии, отвечающие работе машины, образуют некую структуру. Дальнейшее совершенствование и упрощение конструкций приводит к апериодическим мозаикам, ставшим классическими (плитки Амманна, мозаики Пенроуза и др.).

В $1961$ году Хао Вангом \cite{Wang} были рассмотрены квадратные плитки с разноцветными сторонами. Разные плитки можно прикладывать
друг к другу сторонами одного цвета. Был поставлен вопрос, существуют ли конечные наборы таких плиток,
с помощью которых могут быть получены только непериодические замощения плоскости. Первым такой набор был построен
Робертом Бергером \cite{Berger}, основная идея состояла в том, что замощения моделировали работу машины Тьюринга. При этом
использовалось несколько тысяч плиток. Позднее были придуманы наборы, содержащие небольшое количество плиток. Например,
интересна конструкция Рафаэля Робинсона \cite{Robinson}. В этой связи можно также упомянуть мозаики Пенроуза и Амманна, состоящие
из многоугольников, которыми, при заданных условиях,  можно замостить плоскость лишь непериодически.
Различные примеры апериодических мозаик есть также в работе \cite{FerniqueIvanovBlvMitrafnv}.

%Впервые такой набор был построен Робертом Бергером \cite{Berger}. В
%дальнейшем были построены более простые примеры, например, Рафаэлем Робинсоном \cite{Robinson}. Широко известна также мозаика Пенроуза.

\begin{figure}[hbtp]
\centering
\includegraphics[width=0.4\textwidth]{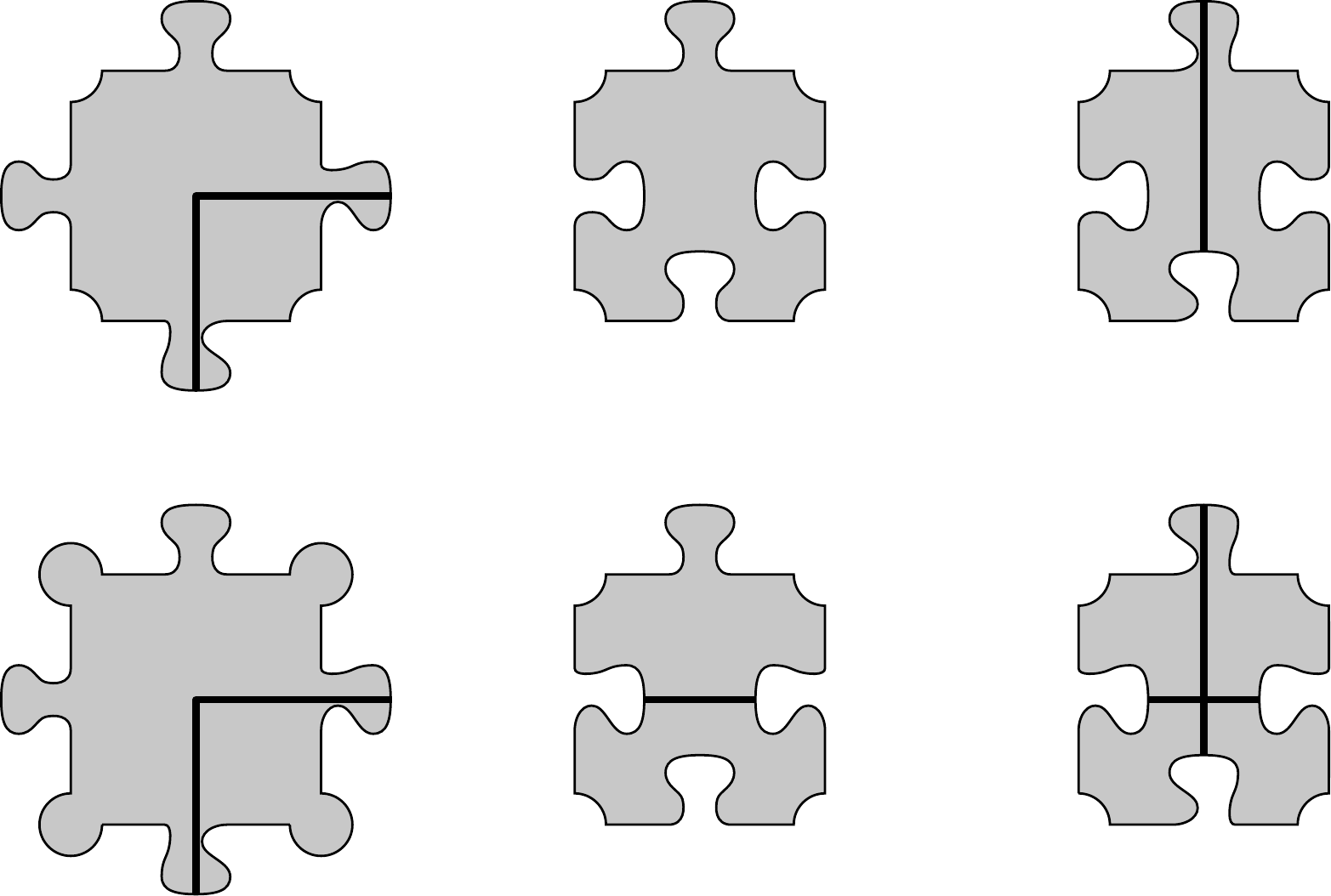}
\caption{Набор плиток Робинсона.}
\label{fig:robinson_bw}
\end{figure}

\subsubsection{Иерархия и апериодичность} \label{dol}

Рассмотрим другой способ получения непериодических замощений. Пусть имеется конечное число типов плиток и мы задаем
правила, по которым из нескольких маленьких плиток можно составлять большие макроплитки тех же типов.

\begin{quote}

{\bf Пример.} {\it Плитки могут быть квадратами $A$ и $B$, при этом, чтобы составить квадрат $A$ второго уровня, нужно взять
четыре квадрата $A$, $A$, $A$, $B$. А чтобы составить квадрат $B$, нужно взять четыре квадрата $B$ $A$ $B$ $A$.}

\end{quote}

Таким образом, получается иерархическая система. Каждую плитку можно разбить на требуемое число уровней иерархии.
Можно показать, что получаемое замощение будет непериодично. Аналогичный способ построения используется в подстановочных
системах, например, c помощью подстановок $1\rightarrow 10$, $0\rightarrow 01$ получается бескубное слово
$$1001011001101001 \dots$$

\medskip

Оказывается, язык граничных условий и язык иерархий схожи. А именно, любую иерархическую систему можно задать конечным числом граничных условий.

Пусть имеется конечное число типов плиток, причем заданы правила иерархии, по которым плитка уровня $n$ составляется из
нескольких плиток уровня $n-1$. Тогда для начального набора плиток первого уровня можно задать конечную систему
граничных условий так, чтобы задавалась мозаика, получаемая при иерархическом способе задания.

Иерархичность системы плиток гарантирует непериодичность замощения. Вследствие этого можно задавать с помощью граничных
условий мозаики, которые будут с гарантией непериодическими.

\subsubsection{Демо-версия доказательства}

\begin{figure}[hbtp]
\centering
\includegraphics[width=0.8\textwidth]{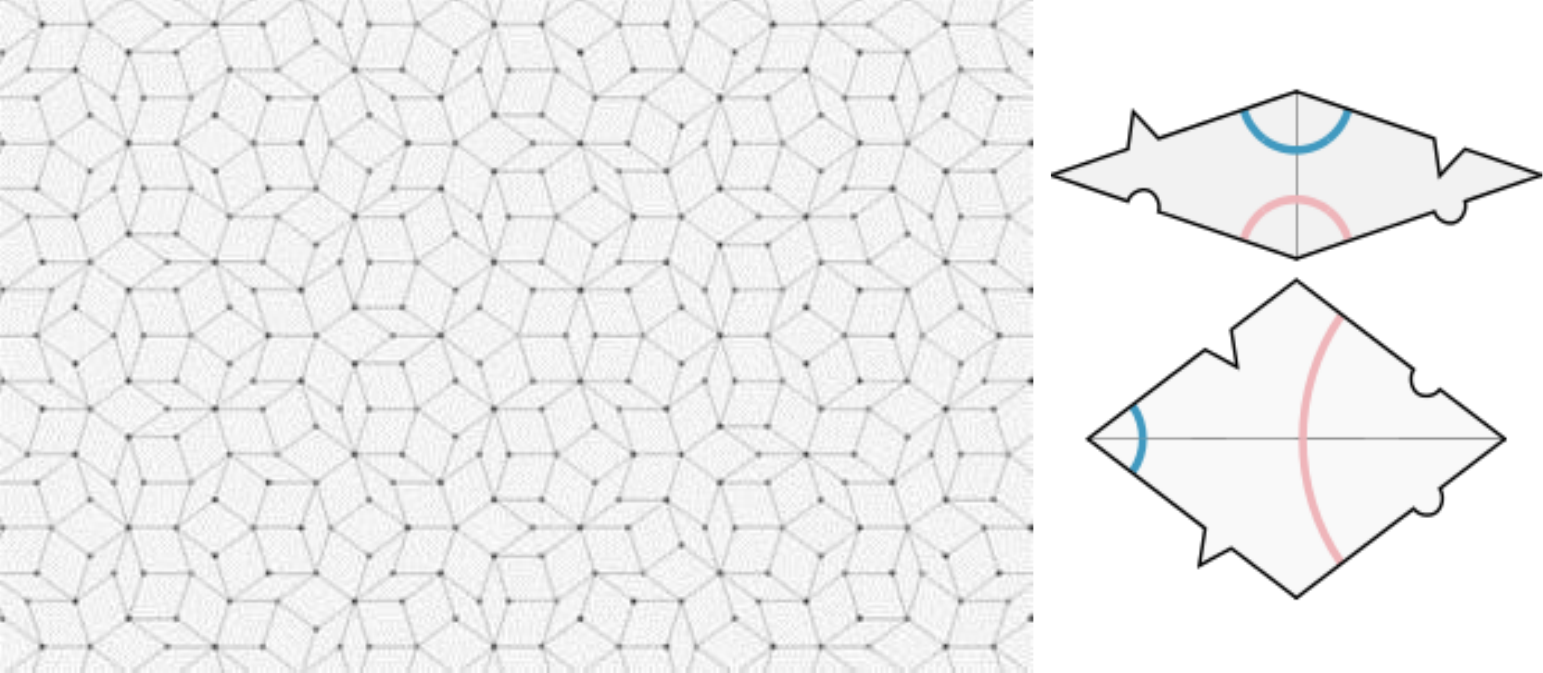}
\caption{Мозаика Пенроуза.}
\label{penrose}
\end{figure}

Рассмотрим одну из классических мозаик -- знаменитую {\it мозаику Пенроуза} (см. рисунок~\ref{penrose}).
\medskip

\paragraph{\bf Конструкция мозаики Пенроуза.}  Используются плитки двух видов -- толстый и тонкий ромбы. Есть граничные условия: стороны каждого ромба раскрашены в две пары цветов. Соприкасаться два ромба могут только сторонами двух цветов, образующих пару. На рисунке~\ref{penrose} цвета обозначены внешними и внутренними насечками двух типов.
\medskip

\medskip

\paragraph{\bf Конструкция полугруппы мозаики Пенроуза.}
Легко видеть, возможно конечное количество типов узлов-вершин где
сходятся несколько плиток. Кроме того, ребра в мозаике могут иметь десять возможных направлений. Можно выписать все возможные типы узлов и обозначить
их буквами алфавита. Теперь последовательность букв (слово) будет кодировать
последовательность узлов, которые мы проходим вдоль пути.

Для каждого узла введем также два параметра: по какому ребру мы в него вошли и по какому ребру мы вышли (включая информацию о цветах ребер). Теперь можно расширить алфавит, добавив буквы для всевозможных сочетаний параметров для разных типов узлов. Некоторые последовательности букв с гарантией не смогут представлять путь на мозаике (например, если ребро, по которому мы пришли в узел не соответствует по цвету ребру, по которому мы вышли из предыдущего узла). Такие последовательности мы будем заносить в список запрещенных.

Некоторые пути на мозаике можно объявить {\it эквивалентными}: например путь по двум соседним ребрам ромба эквивалентен пути по другой паре сторон. Можно выписать
все такие эквивалентности для разных вариантов получающихся узлов и получить
полный конечный список. Тогда, с помощью локальных замен, можно переводить одни пути в другие.

Таким образом, с мозаикой Пенроуза можно связать конечно определенную полугруппу, где словам соответствуют пути. Запрещенные пути -- это нулевые слова. Можно также запретить пути, которые не являются кратчайшими между двумя узлами. Для этого нужно сначала внести в список запрещенных короткие такие пути. После этого можно показать что с помощью локальных замен любой некратчайший путь приводится к виду содержащему запрещенный короткий подпуть.

В случае, если бы можно было бы показать, что любой путь, не вкладывающийся в мозаику, приводится к нулю с помощью указанных локальных правил, получающаяся полугруппа была бы нильполугруппой, так как на мозаике не лежит периодических путей.

\medskip

\paragraph{\bf Почему нужна другая мозаика.}
Проблема заключается в том, что на мозаике Пенроуза есть пути, которые недостаточно сильно меняются локальными заменами, то есть, недостаточно ``извиваются''. Это приводит к тому, что можно сконструировать путь, каждый локальный кусок которого может быть вложен в мозаику, но весь путь не может быть вложен. Локальные замены меняют его незначительно и преобразовать его в достаточной мере, чтобы диагностировать несоответствие мозаике, не получается.

В целом, каждый путь можно трактовать как массив информации о его окрестности. Когда мы производим локальные замены, происходит перенос информации вдоль пути. В случае, если пути можно шевелить незначительно, канал переноса информации будет ограничен, что не позволит выявить ситуацию, когда длинный кусок пути не может являться частью мозаики. (Например, когда путь -- это степень разрешенного слова.)

\subsection{Язык контактов vs язык соотношений. Подклейки.}
Итак, если осуществлять нашу программу на базе классических мозаик, то возникают трудности, связанные с тем, что в некоторых направлениях геодезические пути не поддаются изгибу и вокруг них ничего из переброски не наращивается. А именно достаточно протяженный в двух измерениях кусок мозаики обеспечивает вычислительный процесс.

Перевод с языка соотношений (т.е. движения геодезического пути) на язык контактов вещь более сложная и обеспечить его возможность не так просто. Чтобы локальное шевеление геодезической заполнило достаточную область, пространство должно быть {\it равномерно-эллиптическим}. Пространство называется {\it равномерно-эллиптическим}, если любые две точки $A$ и $B$ на расстоянии $D$ соединяются системой геодезических,  образующих диск ширины $\lambda\cdot D$ для некоторой глобальной константы $\lambda>0$. (Априори не очевидно, что такое пространство существует.) Такое пространство следует собрать из мозаики и установить аналог теоремы Гудмана-Штраусса для этой сборки.

Равномерно-эллиптическое пространство строится из подстановочной системы, указанной на рис. \ref{fig:rule} при этом следует озаботиться тем, чтобы степени узлов были бы ограничены, что усложняет конструкцию, вынуждая сделать композицию разных подстановок.

%Далее, когда путь проходит через узел и по обе стороны от него есть пара плиток более низкой ступени иерархии, сходящиеся в этом узле, необходимы {\it подклейки}. Бесконечное число иерархических уровней приводит к тому, что к ребру может примыкать неограниченное число подклеек. При этом необходимо избежать следующей ситуации: мы вышли в подклейку и вернулись. Далее продеформировали участок пути в подклейку так, что его ``голова'' снова коснулась исходного ребра, а вход и выход в ``подклейку'' не изчезли. Далее возникший участок по ребру может уйти в иную подклейку там снова изогнуться и коснуться исходного ребра и т.д. Конструкция специально построена таким образом, чтобы этого избежать.
%То есть {\em все организовано так, что путь может вернуться из подклейки ТОЛЬКО в исходную точку}. Далее размеры подклеек к плиткам, подклеек к подклейкам и т.д. экспоненциально убывают  достаточно быстро. Тем самым, помимо всего прочего, обеспечивается гомотопическая тривиальность комплекса, образованного ``подклейками''.

Мы хотим добиться того, что никакой узел $X$ на пути из A в B не являлся обязательным пунктом посещения (узким местом). Для этого вводятся  {\it подклейки}. Для некоторых путей $AXB$ мы будем вклеивать макроплитку  $AXBZ$, где $Z$ -- новая создаваемая вершина, вне плоскости $AXB$. На дальнейших этапах разбиения проводятся подклейки к подклейкам. Если путь устроен так, что мы входим внутрь подклейки, а потом возвращаемся из нее, то участок на пути, находящийся вне базовой плоскости все время является односвязным.

\begin{quote}
Это, в частности, означает, что нельзя деформировать участок пути в подклейку так, чтобы он частично стал лежать в базовой плоскости. Кроме того, размеры подклеек к плиткам, подклеек к подклейкам и т.д. экспоненциально убывают достаточно быстро. Тем самым, помимо всего прочего, обеспечивается гомотопическая тривиальность комплекса, образованного ``подклейками''.

\end{quote}

\subsection{Последовательная канонизация. Завершение доказательства.}
Доказательство завершается так. Рассмотрим слово. Надо его либо привести к кодировке, соответствующей пути на комплексе (тем самым оно непериодично), либо к нулю.
Оно последовательно приводится к {\it $k$-каноническому виду} и при этом $k$ растет, т.е. оно состоит из участков границ плиток $k$-го уровня иерархии, кроме начала и конца, целиком содержащихся в плитке уровня $k-1$ (возможно в подклейке $(k-1)$-го уровня) и при этом являющихся $(k-1)$-каноническими.

Далее с помощью соотношений (в том числе отвечающих выходу в подклейку) проверяется согласование соседних участков пути, расположенных на границах плиток $k$-го уровня, их нахождение в $(k+1)$-ом уровне или возможность находиться в двух соседних плитках $(k+1)$-го уровня. Отсутствие согласования означает, что пути с такой кодировкой не существует на комплексе -- т.е. в этом случае слово можно привести к нулю.

В начале рассматривается плоский процесс. Если имеются команды входа в подклейку и выхода из нее, то мы преобразуем слово между ними и производим сокращение -- убирается один выход в подклейку.

Согласование означает возможность преобразования к $(k+1)$-каноническому виду, после чего процесс повторяется. Он заканчивается канонической формой слова. В нашем случае каноничность не означает однозначности, поскольку можно выбирать разные стороны плиток. Наверное, существует такой выбор -- но он более сложен и ситуация похожа на теорию канонической формы элементов в гиперболических (или в более общем случае, в группах с неположительной кривизной) развитую И.~А.~Рипсом, когда вначале строится предканоническая форма (для данного уровня иерархии), потом осуществляется выбор.

\subsection{Специфика геометрических методов в полугрупповом случае}
Интерпретация соотношений в полугруппах через ``rewriting diagram'' похожа на применение диаграмм ван Кампена в теории групп.
Рассмотрение мозаик также родственно работе с диаграммами ван Кампена в теории групп.
Использование путей по ребрам в плиточных замощениях и групп связанных с ними, встречается в работах Дж.~Конвея \cite{Conway}.
С помощью групповой техники он показал отсутствие различных разбиений. Вот модельный

\medskip
{\bf Пример.} {\it Рассмотрим шахматную доску, из которой вырезаны две противоположные угловые клетки. Тогда ее нельзя разбить на доминошки $1\times 2$}.
\medskip

Предположим, что такое разбиение существует. Рассмотрим доминошки как клетки в диаграмме ван Кампена. Букве $a$ сопоставим вертикальную стрелку, букве $b$~-- горизонтальную, обратным буквам -- обратные стрелки. Вертикальная доминошка отвечает соотношению $a^2ba^{-2}b^{-1}=e$, а горизонтальная -- $b^2a^{-1}b^{-2}a=1$. Границе квадрата с вырезанными углами отвечает путь $a^7bab^7a^{-7}b^{-1}a^{-1}b^{-7}$. Будем интерпретировать групповые слова как элементы группы $S_3$. Если сопоставить элементу $a$ транспозицию элементов $1$ и $2$, а элементу $b$ -- транспозицию $2$ и $3$, то $a^2=b^2=e$ и соотношения выполнены, в то же время как $a^7bab^7a^{-7} b^{-1}a^{-1}b^{-7}=(ab)^4=ab\ne e$ есть трехчленный цикл.

Тем же методом решаются многие другие олимпиадные задачи на раскраску, но в то же время и такие, которые раскраской не делаются (подробности см. \cite{BelvIvanvMalstvMitrfnvKharitnv}).

Тем не менее мозаики и диаграммы ван Кампена существенно различаются. Например, свойства величин углов имеют далеко не полное отражение. Аналогия между мозаиками и группами довольно глубокая, но не вполне формализована. Имеется классический результат об алгоритмической разрешимости проблемы равенства слов в группе с одним соотношением. Имеется и аналогичный вопрос об алгоритмической разрешимости проблемы разбиения плоскости транслятами одной фигуры. Для связных фигур такая алгоритмическая разрешимость доказана, а в общем случае вопрос остается открытым.

Далее, при разбиении квадрата на домино возникают подквадратики $2\times 2$, разбитые на пары доминошек. Такую пару можно повернуть на $90^o$ и сделать {\it флип}. Такими флипами можно от одного разбиения перейти к любому другому. Аналогичный факт верен и для разбиения на $k$-миношки и для многих других разбиений. Данный факт родственен конечной порождаемости группы $\pi_2$ комплекса (с учетом действия фундаментальной группы $\pi_1$). Хорошо бы прояснить аналогию.

В полугруппах геометрические идеи работают не совсем так, как в группах, и в нашем случае эффекты ``углов'' лучше представлены.   По всей видимости, полугрупповая теория может пролить свет на мозаики. Более того, есть некая близость теории полугрупп к кольцам как к ``квантовой взвеси'', так что и кольцевая теория, если будет построена, даст дополнительное понимание.
Возможно, наш  подход окажется полезным и для других построений в полугруппах и кольцах.

\subsection{Апериодические замощения в нашей конструкции}
Как уже говорилось,  существуют конечные наборы многоугольников (плиток) которыми можно
замостить плоскость лишь непериодически. Методы построения таких мозаик, как правило, опираются на иерархические правила построения: задаются универсальные правила построения плиток уровня $n+1$ из плиток уровня $n$, для нескольких типов плиток $A_1,\dots,A_k$.
Х.~Гудман-Штраусс \cite{GoodmanStrauss}  показал, что иерархические мозаики можно получить, задав конечное число локальных правил.
Таким образом, локальные условия (конечность набора) могут приводить
к глобальному эффекту (непериодичности замощения).

Итак, основной задачей является конструирование мозаики, в которой указанные выше трудности были бы решены. Для этого она должна обладать несколькими свойствами:

\medskip
\begin{description}
  \item[1. Локальная конечность] Речь идет о конечности числа возможных типов узлов, конечности изначально задаваемых пар эквивалентных путей, а также конечности списка запрещенных путей.
  \item [2. Возможность ``шевеления'' любого пути] Любой достаточно длинный путь, соединяющий узлы $A$ и $B$ может быть переведен локальными заменами в другой путь, отличающийся достаточно сильно от начального.

      Более точно, пусть длина пути равна $n$ и точка $M$ -- его середина (или такая точка, что расстояния вдоль путей $AM$ и $MB$ различаются на 1). Пусть $AM'B$ -- эквивалентный путь, $M'$ -- соответствующая $M$ точка. Пусть $R_{AMB}(n)$ -- максимальное расстояние от $M$ до $M'$ по всем эквивалентным путям $AM'B$.
      Требуемое нами условие означает, что $R_{AMB}(n)$ является неограниченной возрастающей функцией от $n$.

  \item [3. Апериодичность] На мозаике не должно быть путей, отвечающих периодическим словам.
\end{description}

 \medskip
Как уже говорилось, мы используем геометрическую интерпретацию для алгебраических построений. Запрет для
двух (или более) плиток находиться рядом друг с другом схож с запретом для двух букв стоять рядом в разрешенном слове.
Возникает интерпретация слова как последовательности плиток на выложенной мозаике.
Апериодичность мозаики приводит к непериодичному характеру таких ``плиточных слов''.  В свою очередь,
непериодических замощений можно добиться, если применять иерархический способ построения.

В связи с этим используются  языки плиточных примыканий и иерархий. Эти языки во многом схожи, например, заданные
правила иерархии, когда плитки $A_1,\dots,A_k$ уровня $n+1$ составляются из плиток $A_1,\dots,A_k$ уровня $n$, можно выразить
с помощью конечного множества локальных правил для набора $A_1,\dots,A_k$. Эти локальные правила будут порождать те же непериодические мозаики что и исходные правила иерархии.

\medskip

Дальнейшее развитие этих языков приводит к появлению более универсального языка путей на графе. То есть плиточная мозаика
рассматривается как граф, где вершины это узлы мозаики, а ребра -- границы плиток. В этом смысле понятие плитки
можно обобщить, рассматривая их уже как локальные подграфы из которых, с помощью локальных правил, можно составлять граф,
покрывающий плоскость.

Аналогом буквы будет тип вершины графа, аналогом слова -- путь,
проходящий через несколько вершин. Аналогом соотношения будет эквивалентность между путями с общими концами: например,
в простом $4$-цикле $ABCD$ выполнено соотношение $ABC=ADC$. Помимо таких, есть также мономиальные соотношения, выражающие
идею о невозможности существования какого-то пути на мозаике. Также, для обеспечения необходимого контроля над множеством
ненулевых слов, вводятся мономиальные соотношения, обнуляющие слова, соответствующие некратчайшим путям.
Оказывается, что при этом можно обойтись конечным числом соотношений. Немономиальные соотношения, при этом, не меняют длины слова.

\medskip

Язык путей на мозаике-графе позволяет выразить те же концепции и определить те же мозаики, что и языки иерархических плиток или граничных условий.
Таким образом возникает связанная с мозаикой конечно определенная полугруппа с набором свойств, индуцированных мозаикой.

\begin{figure}[hbtp]
\centering
\includegraphics[width=0.3\textwidth]{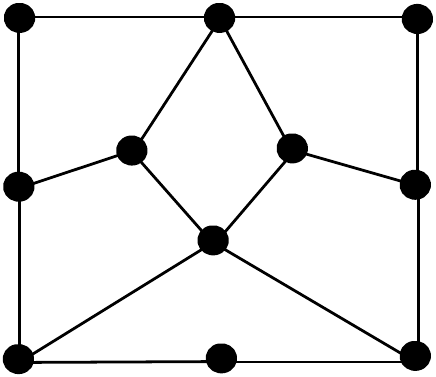}
\caption{Правило разбиения.}
\label{fig:rule}
\end{figure}

Для построения нильполугруппы используется мозаика сгенерированная с помощью следующего иерархического правила разбиения (рисунок~\ref{fig:rule}).

От мозаики нужно потребовать ряд дополнительных свойств. В частности, любой длинный путь должен допускать
возможность ``шевеления'' то есть, локальных преобразований над ним, позволяющих в достаточной мере менять его. Для
достижения этого свойства к плоской мозаике производятся ``подклейки'', представляющие собой небольшие плоские подграфы, не лежащие
в исходной плоскости, и позволяющие обходить ``узкие места'' на исходном плоском графе. Структура подклеек также имеет иерархическую
природу и так же может быть задана на языке преобразований путей конечным образом.

В итоге мы получаем граф, обладающий набором важных для нас свойств, в котором любой путь, соединяющий произвольные точки $A$ и $B$,
является членом семейства геодезических эквивалентных друг другу путей, соединяющих эти точки. Причем эквивалентность двух путей
из этого семейства может быть получена путем цепочки локальных преобразований переводящих один путь в другой. ``Шевеление''
пути играет роль передачи информации. Фактически, определяющие соотношения задают правила передачи информации по пути.
Если задан произвольный длинный путь, мы можем начать работать над ним, совершая локальные преобразования. При этом возможна альтернатива:

\begin{enumerate}
  \item В результате этой работы можно получить внутри некратчайший подпуть, либо подпуть указанный в числе запрещенных.
В этом случае наш путь представляет нулевое слово.
  \item В результате этой работы мы восстанавливаем некоторый кусок мозаики, внутри которого лежит пучок геодезических путей, эквивалентных нашему.
\end{enumerate}

\medskip

Мозаика не может содержать в себе путей, выражаемых периодическим словом. То есть, все периодические слова могут быть приведены к нулю
 локальными преобразованиями. При этом геодезические пути, лежащие на мозаике не приводятся к нулю, и могут иметь любую длину.
Таким образом полугруппа, соответствующая построенному графу-мозаике, будет конечно определенной нильполугруппой.

\subsection{Конечно определенные полугруппы} \label{semigroups}
%\medskip
Рассмотрим алфавит $\Omega$, состоящий из конечного числа букв. Слова, составленные из букв образуют полугруппу
относительно операции приписывания одного слова к другому. Кроме того, есть специальная буква $0$ (ноль), такая что
$W0=0W=0$ для любого слова $W$ из полугруппы. Определяющим соотношением в полугруппе считается равенство вида
$W_1=W_2$, где $W_1$ и $W_2$ -- некоторые слова, причем одно из них может быть нулем (или нулевым словом). В этом случае,
соотношение называется {\it мономиальным}. Для натурального $n$ степенью слова $W^n$ называется слово $WW\cdots W$, где $W$ выписано
$n$ раз подряд. Элемент полугруппы (слово) называется {\it нильэлементом}, если существует натуральное $n$, такое что
$W^n=0$. Если все элементы полугруппы являются нильэлементами, то вся полугруппа называется нильполугруппой.
Разумеется, наша полугруппа содержит {\it ноль}.

Мы пользуемся геометрической интерпретацией: {\it буквам} отвечают вершины специального
графа, а {\it словам} -- пути в этом графе. И если слово не может быть представлено путем на графе, оно всегда будет
приводиться к нулю.

%\medskip

\subsubsection{Плитки и непериодические замощения} \label{plitki}

Начнем с чистой геометрии.
Мы задаем конечный набор граничных условий (как можно прикладывать плитки друг к другу),
и через задание локальных условий достигается глобальный эффект. Задание граничных условий схоже с заданием определяющих
соотношений в полугруппе, этим объясняется интерес к плиткам и замощениям.

\subsubsection{Узлы на мозаике и пути по границам плиток} \label{ways}

Помимо языков иерархических систем и граничных условий есть еще один подход. Можно рассматривать мозаики с
точки зрения путей по границам плиток. Это естественный шаг для перехода к полугруппе, так как пути по мозаике
логично отождествляются со словами в подходящем алфавите.

Назовем {\it узлами} вершины мозаики, где сходится несколько плиток. Ясно, что в плоской мозаике с конечным набором
типов плиток возможно конечное число видов узлов. Обозначим их буквами конечного алфавита. Последовательность
букв (слово) соответствует пути на мозаике. Некоторые слова могут вообще не встречаться на мозаике, а другие встречаться в разных местах.

Третий подход к конструированию мозаик состоит в определении конечного списка невозможных путей, а также задании конечного
числа пар эквивалентных путей. С помощью такого языка можно задавать мозаики так же, как и с использованием
иерархических систем плиток или граничных условий на плитки.

Мы будем строить такую мозаику в два приема. На первом этапе построим ее {\it плоскую часть}, которую будем считать базовой. Она будет строиться в виде иерархического графа. На втором этапе мы покажем, как к построенному графу производятся {\it подклейки} -- плоские и конечные подграфы, лежащие в другой плоскости, но имеющие с базовым графом несколько общих ребер. В результате получится 2-комплекс вложенный в трехмерное пространство.

\medskip

Подклейки нужны для того, чтобы появился маршрут, альтернативный проходу по ребру, к которому производится подклейка. В этом случае, становится возможным перенос информации вдоль пути и мы можем выявить ситуации с невозможным на мозаике путем.

\medskip

В дальнейших главах мы построим мозаику, где выполнены указанные выше свойства.

%\medskip

\subsection{Схема зависимости материала и состав глав}

Доказательство основного результата достаточно объемное. В этом параграфе мы кратко опишем состав глав, составляющих работу.

\medskip

Во {\bf введении} обсуждаются проблемы бернсайдовского типа и конечной определенности и связанные с этим вопросы.

В главе~\ref{ScFnDefInfGen} обсуждается общая стратегия доказательства а также смежные вопросы.

В главе~\ref{geom} (``Геометрическая структура комплекса'') обсуждаются чисто геометрические свойства комплекса, лежащего в основе построения.

В последующих главах продолжается построение конечно определенной нильполугруппы. 
Вводится параметризация, благодаря которой структура построенного геометрического комплекса оказывается связанной со структурой полугруппы. В частности, вводится алфавит букв, кодирующих вершины и ребра на комплексе. Детерминированность комплекса позволяет ввести определяющие соотношения на коротких отрезках путей. Кратчайшим путям на комплексе будут соответствовать ненулевые слова во введенной полугруппе.

Далее с помощью геометрических и комбинаторных свойств комплекса доказывается, что произвольное слово содержащее девятую степень приводится к нулю с помощью определяющих соотношений.

\medskip

%В настоящей статье развивается новый подход к контролю над вводимыми соотношениями. Он основан на
%следующих соображениях.

%\subsection{Плитки, мозаики и полугруппы, уточнение плана доказательства}

\section{Геометрическая структура комплекса} \label{geom}

\subsection{Иерархическое построение}

Построим {\it комплекс уровня $n$} с помощью итерационного процесса. На каждом шаге мы будем иметь дело с графом, причем простые циклические пути из четырех ребер будем называть {\it плитками}.

\medskip

\begin{figure}[hbtp]
\centering
\includegraphics[width=0.4\textwidth]{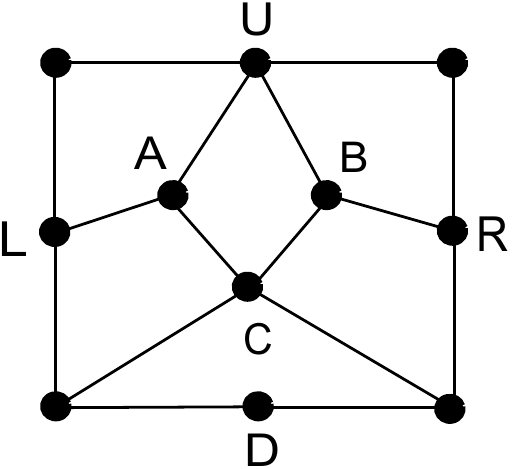}
\caption{Комплекс второго уровня.}
\label{fig:firststep}
\end{figure}

{\it Комплекс уровня $1$} представляет собой простой цикл из четырех вершин $CUL$, $CUR$, $CDR$, $CDL$. {\it Комплекс уровня $2$} это граф на рисунке~\ref{fig:firststep}. 
%Вершины $A$, $B$, $C$ будем называть {\it внутренними}, а вершины $U$, $L$, $R$, $D$ -- {\it боковыми}. 

То есть,  на втором этапе у нас имеется шесть плиток. Будем называть их, согласно их положению, левой верхней, левой нижней, правой нижней, правой верхней, средней, нижней.

Будем считать, что вершины $A$, $B$, $C$, $U$, $L$, $R$, $D$ имеют {\it глубину} 0, изначальные четыре вершины (углы) имеют глубину равную $-1$. Будем также называть их {\it угловыми}. Кроме этого, угловые вершины будут появляться при операции подклейки (см ниже).

\medskip

С комплексом будут производиться итерации двух типов: {\it разбиение} и {\it подклейка}.

\begin{figure}[hbtp]
\centering
\includegraphics[width=0.5\textwidth]{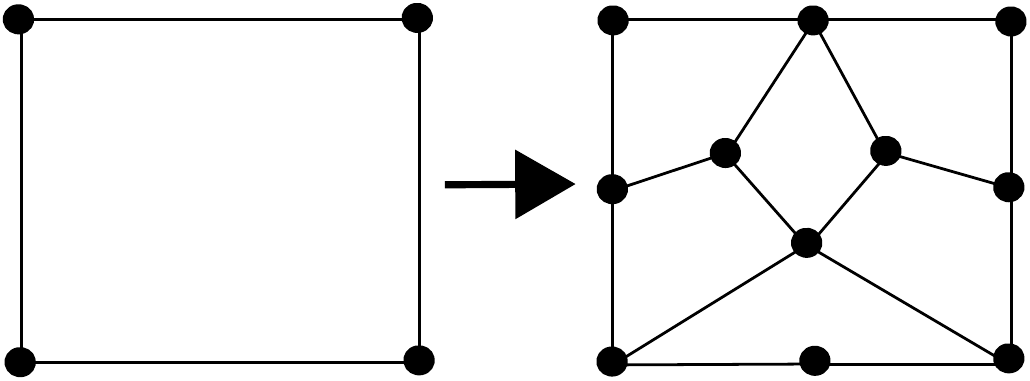}
\caption{Разбиение.}
\label{fig:subst}
\end{figure}

\begin{figure}[hbtp]
\centering
\includegraphics[width=0.5\textwidth]{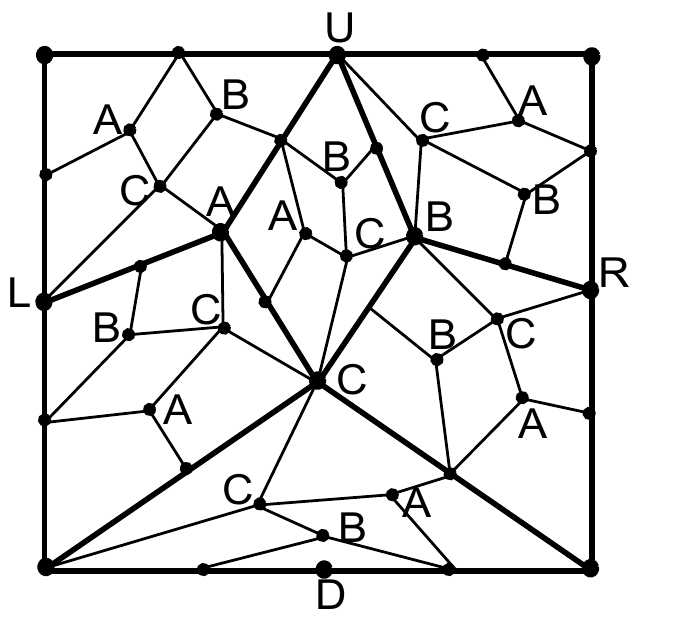}
\caption{Макроплитка третьего уровня. Отмечены типы внутренних вершин ($A$, $B$, $C$)}
\label{fig:level2}
\end{figure}

\begin{definition}
{\it Макроплиткой} уровня $n$ назовем: для $n=1$ -- обычную плитку (простой $4$-цикл), которую будем также называть {\it минимальной}; для $n>1$ -- результат применения следующей операции {\it разбиения} к макроплитке уровня $n-1$: макроплитка разбивается на шесть плиток согласно правилу на рисунке~\ref{fig:subst}. создаются новые четыре вершины на границе разбиваемой плитки (в серединах ее сторон) и три вершины внутри нее. Одна из четырех возможных ориентаций определяется согласно положению (одному из шести) разбиваемой плитки в ее родительской макроплитке (рисунок~\ref{fig:level2}).

Одна из четырех возможных ориентаций определяется согласно положению (одному из шести) разбиваемой плитки в ее родительской макроплитке второго уровня (рисунок~\ref{fig:level2}).

Проводимые при разбиении ребра будем называть {\it внутренними ребрами, принадлежащими} разбиваемой макроплитке. {\it Уровнем ребра} будем называть уровень макроплитки, которой оно принадлежит. То есть, принадлежность ребра определяется раз и навсегда, при его появлении. {\it Типами} внутренних ребер будем называть $8$ видов ребер образующихся при разбиениях (рисунок~\ref{fig:inneredges}). Также для каждого внутреннего ребра определим две стороны $A$ и $B$, расположенные как указано на рисунке. Есть также $4$ типа {\it граничных} ребер (это ребра изначального цикла из $4$ вершин), левое, правое, верхнее и нижнее. Пусть они также пронумерованы от $9$ до $12$. Граничные ребра комплекса относятся к этим типам. Также к ним относятся границы создаваемых подклееных макроплиток (см ниже). Будем считать, что граничные ребра принадлежат макроплиткам, сторонами которых они являются.

{\bf Замечание.} Далее в тексте термин {\it ребро} означает ребро макроплитки, а не ребро графа, если только из контекста не следует обратное. В то же время,  {\it входящее } или {\it выходящее} ребро -- как раз ребро графа, входящее в данную вершину.

\begin{figure}[hbtp]
\centering
\includegraphics[width=0.4\textwidth]{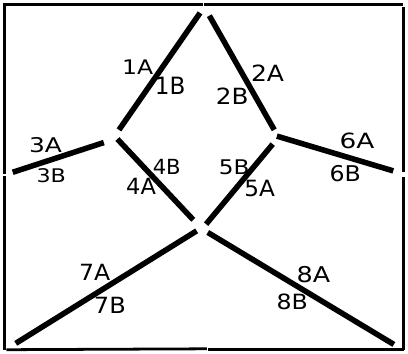}
\caption{Восемь типов внутренних ребер c ориентацией.}
\label{fig:inneredges}
\end{figure}

\end{definition}

Определим {\it типы вершин} на комплексе.
\begin{definition}
Все вершины, встречающиеся на комплексе, мы разделим на следующие категории:

1) {\it Угловые}. (Лежащие в углах подклееных макроплиток или всего комплекса). {\it Тип угловой вершины } определим как один из четырех вариантов углов, в котором она может находиться: $\mathbb{CUL}$, $\mathbb{CUR}$, $\mathbb{CDR}$, $\mathbb{CDL}$. (Corner Up-Left и так далее.) Угловые вершины {\it принадлежат} макроплитке, где они являются углами.

2) {\it Краевые}. (Лежащие на стороне подклееной макроплитки или всего комплекса). Каждая такая вершина лежит в середине стороны некоторой макроплитки, прилегающей к краю. Тип краевой вершины  определим в соответствии с тем, серединой какой стороны в этой макроплитке она является: $\mathbb{L}$, $\mathbb{R}$, $\mathbb{D}$, $\mathbb{U}$.
Краевые вершины {\it принадлежат} макроплитке, на краю которой они лежат.

3)  {\it Внутренние}. В этой категории определим три типа внутренних вершин: $\mathbb{A}$, $\mathbb{B}$, $\mathbb{C}$, отвечающих внутренним узлам макроплиток, эти вершины создаются {\it внутри} разбиваемой макроплитки. Будем считать, что эти вершины {\it принадлежат} данной макроплитке.

4)  {\it Боковые}. (Лежащие на границе между двумя макроплитками, на внутреннем ребре).
  Типы боковых вершин будут соответствовать всем упорядоченным парам из множества $\{ \mathbb{U},\mathbb{R},\mathbb{D},\mathbb{L}\}$. А именно: $\mathbb{DR}$, $\mathbb{RD}$, $\mathbb{DL}$, и так далее. Будем считать, что в упорядоченной паре первой называется буква, соответствующая $A$-стороне внутренего ребра, на котором лежит вершина.  Все боковые вершины создаются в середине стороны разбиваемой макроплитки. 
Тип боковой вершины определяет, серединой каких именно сторон она является в двух макроплитках, где она является серединой сторон. Это как раз те макроплитки, которые разбивались при создании данной вершины.
Будем считать, что боковые вершины {\it принадлежат} макроплитке, которой принадлежит данное внутреннее ребро.

В дальнейшем, вершину типа $\mathbb{A}$ будем, для простоты, называть $A$-вершиной или $A$-узлом. Аналогично для других типов. Вообще, вершины на графе мы также будем называть узлами.

Заметим, что тип вершины определяется при ее создании и в дальнейшем не меняется. Например, если вершина  была боковой вершиной и после какой-то подклейки оказалась на краю подклееной макроплитки, она все равно останется боковой и принадлежащей той макроплитке, внутреннему ребру которой она принадлежит.
\end{definition}

\medskip

Определим {\it глубину} новых вершин. Все создаваемые при разбиении вершины получают глубину на 1 больше, чем максимальная глубина вершины до разбиения.
Из построения ясно, что у любого ребра макроплитки хотя бы один из концов является вершиной, созданной на предыдущем шаге. Следовательно, создаваемая вершина в середине ребра получает глубину на $1$ большую, чем, максимальная глубина двух концов ребра.

\begin{definition}
Операция {\it разбиения}, примененная ко всему комплексу есть результат одновременного применения разбиения ко всем плиткам первого уровня в комплексе.
\end{definition}

\begin{definition}
Операция {\it подклейки}. При {\it подклейке} мы рассматриваем все пути $X_1X_2YZ_2Z_1$ из пяти  вершин (и четырех ребер), такие что:

1) Вершины $X_1$, $Y$, $Z_1$ {\bf не} являются тремя углами из четырех никакой макроплитки;

2) Вершины $X_1$, $Z_1$ являются боковыми или краевыми вершинами глубины $k-1$, где $k$ -- максимальная глубина вершины в комплексе;

3) Вершина $X_2$ является серединой $X_1Y$, а $Z_2$ -- серединой $Z_1Y$, то есть,  $X_2$ и $Z_2$  являются боковыми или краевыми вершинами глубины $k$, созданными в комплексе самыми последними, при последней операции {\it Разбиения}.

4) Вершина $Y$ имеет глубину $k-2$.

5) Путь $X_1X_2YZ_2Z_1$ выбран так, что уровень ребра, на котором лежит вершина $X_1$ больше уровня ребра вершины $Z_1$.

Допустим, ребра, на которых лежат $X_1$ и $Z_1$, имеют один уровень. В этом случае мы считаем, что ребро $X_1$ имеет меньший номер в 
нумерации входящих ребер, которая определяется индуктивно по глубине, см ниже.

%Учитывая, что в макроплитке только одно ребро каждого типа, возможны следующие случаи:

%случай А) ребра, на которых лежат $X_1$ и $Z_1$ принадлежат макроплиткам из разных подклееных областей (одна из областей может быть изначальной базовой макроплиткой). Рассмотрим эти области, как макроплитки. Заметим, что эти макроплитки имеют разный уровень, так как при операции подклейки вершина может попасть только в одну подклеиваемую макроплитку за одну операцию. Выберем за $X_1$ то ребро, макроплитка которого имеет больший уровень (или, другими словами, чья подклееная плоскость появилась раньше).

%случай B) $X_1$ и $Z_1$ в одной плоскости и лежат на ребрах разного типа. Тогда $X_1$ выбирается как вершина с большим номером типа ребра.

%случай C) $X_1$ и $Z_1$ лежат на одном ребре. Тогда для внутреннего ребра $X_1$ выбирается как вершина лежащая ближе к краю макроплитки, если расстояние одинаково, то до середины верхней стороны макроплитки. Для граничного ребра $X_1$ выбирается как предшествующая $Z_1$ при обходе контура по часовой стрелке.

\end{definition}

\medskip

\begin{figure}[hbtp]
\centering
\includegraphics[width=0.5\textwidth]{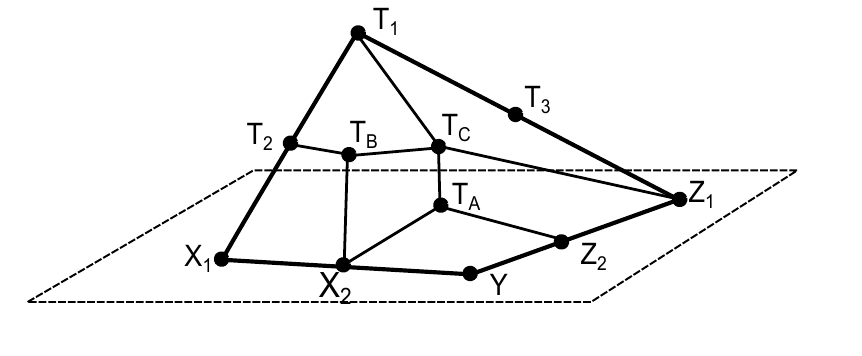}
\caption{Подклейка.}
\label{fig:pasting}
\end{figure}

Далее, для каждого такого пути создаются шесть новых вершин $T_1$, $T_2$, $T_3$, $T_A$,  $T_B$,  $T_C$, не лежащих в плоскости $X_1YZ_1$, а также проводятся новые ребра $X_1 T_2$, $X_2 T_A$, $X_2 T_B$, $T_2 T_B$, $T_C T_B$, $T_C T_A$, $T_2 T_1$, $T_3 T_1$, $T_3 Z_1$, $T_C Z_1$, $T_A Z_2$  и появляется новая {\it подклееная } макроплитка $X_1YZ_1T_1$ (рисунок~\ref{fig:pasting}).

\medskip

Созданные вершины будем называть {\it подклеенными}, причем $T_1$ присваивается глубина $k-1$, а вершинам $T_2$, $T_3$, $T_A$,  $T_B$,  $T_C$ -- глубина $k$.
Вершину $T_1$ будем считать угловой, вершины $T_2$, $T_3$ -- краевыми, а вершины $T_A$,  $T_B$,  $T_C$ -- внутренними.

Типы вершин ($A$,$B$,$C$,$U$, $UL$ и т.п.) и типы ребер не меняются при проведении подклейки.

Вершину $Y$ будем называть {\it ядром подклейки}. Макроплитку, которой принадлежит ядро подклейки, будем называть {\it базовой плоскостью} подклейки.
Макроплитки, которым принадлежат ребра $YX_1$ и $YZ_1$ (это может быть и одна макроплитка) будем называть макроплитками, к которым подклеивается новая 
подклееная область.

При операции подклейки проводятся новые ребра, являющиеся внутренними ребрами в подклееной макроплитке, а также два граничных ребра - нижнее и правое. Рассмотрим те из проведенных новых ребер, которые одним концом лежат на верхнем и левом ребрах подклеиваемой макроплитки, то есть тех ребрах которыми она подклеивается к остальному комплексу. Будем называть их {\it ребрами входа и выхода} в данную подклейку. Если путь проходит по такому ребру, то мы будем говорить что путь вошел в данную подклееную макроплитку, или вышел из нее, в зависимости от направления пути. Заметим, что ребра входа и выхода могут появиться также после разбиения уже подклееной макроплитки.

\medskip

{\bf Примечание 1.} Заметим, что в подклееной макроплитке по построению определяется верхняя, а, следовательно, и остальные стороны: она выглядит так же, как если бы к плитке $X_1YZ_1T_1$ применили разбиение, считая сторону $X_1Y$ верхней. Также можно считать, что в этой макроплитке $T_2$ середина стороны $T_1X_1$, а $T_3$ -- середина стороны $T_1Z_1$.

\medskip

\begin{definition}
 Комплексами $1$, $2$ и $3$ уровня будем называть макроплитки, соответственно, $1$, $2$ и $3$ уровня.

 Комплексом $n$ уровня для $n\geq 4$ становится комплекс $n-1$ уровня, к которому применены {\it Разбиение} и затем {\it Подклейка}.

Иногда мы будем рассматривать {\it плоские части} комплекса.  Так мы будем называть отдельно рассматриваемые макроплитки, полученные разбиениями изначальной плитки, а также отдельно рассматриваемые подклееные макроплитки. {\it Плоским} путем будем называть путь, полностью лежащий в некоторой макроплитке $T$ ( этот путь не содержит ребер, ведущих в подклееные к $T$ макроплитки). 

 \end{definition}

\medskip

%{\bf Замечание 1.} После разбиения, примененного к комплексу первого уровня, образуется граф на рисунке~\ref{fig:firststep}, на котором нет пяти вершин, где можно было бы выполнить подклейку. То есть, подклейки начинают проводиться после второго разбиения.

\medskip

{\bf Замечание 1.} Начиная с четвертого уровня комплекса (и при построении пятого уровня) путь $X_1YZ_1$ из определения операции подклейки может частично лежать в базовой плоскости, а частично в подклееной части (если $Z_1$ появилась при предыдущей подклейке в качестве вершины $T_2$ или $T_3$). На больших уровнях комплекса весь путь из определения подклейки может лежать полностью в подклееных плитках. Таким образом, возникают подклейки к частям подклееных когда-то макроплиток и т.д. То есть подклееная макроплитка подклеивается либо к одной либо к двум другим макроплиткам (к двум в случае, когда два ребра подклейки принадлежат разным макроплиткам). Заметим также, что уровень обоих этих макроплиток не менее чем на $1$ превосходит уровень подклеваемой макроплитки. 

{\bf Замечание 2.} В итоге мы всегда будем рассматривать комплексы конечного уровня. Предельного перехода применяться не будет.

%{\bf Замечание 3.} Ясно, что при проведении разбиений и подклеек угловая вершина остается угловой, боковая--боковой, а внутренняя (относительно базовой плоскости) -- внутренней. Кроме того, можно заметить, что тип вершины ($A$,$B$,$C$,$U$, $UL$ и т.п.) определяется в момент появления этой вершины при разбиении или подклейке и не меняется при дальнейших разбиениях и подклейках.

\medskip

\begin{lemma}[Об ограниченности роста степени вершины] \label{growth_bound}

Для каждой вершины $Z$ существует такое натуральное $N$, что
начиная с уровня макроплитки $N$, степень (число входящих ребер) вершины $Z$ не меняется, то есть она одинакова для макроплиток уровня $N$ и $N+k$ для любого натурального $k$.
\end{lemma}

\medskip

\begin{proof} Рассмотрим некоторую плитку $T$. У нее четыре угла $X$, $Y$, $P$, $Q$ (рисунок~\ref{fig:tileangle}).

\begin{figure}[hbtp]
\centering
\includegraphics[width=0.5\textwidth]{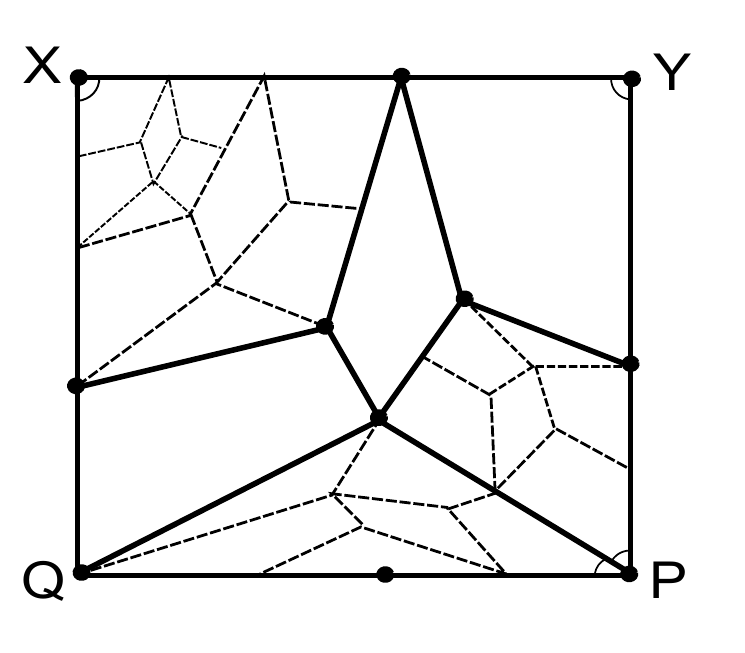}
\caption{Простые углы отмечены дугами.}
\label{fig:tileangle}
\end{figure}

При разбиении некоторые углы разбиваются ребрами, а некоторые нет.

Заметим, что левый и правый верхние углы не будут разбиты ребрами, какой бы ни был уровень макроплитки. Действительно,
левый верхний угол переходит в себя при разбиении, и поэтому никогда не будет разбит ребрами. К правому верхнему углу, на следующем уровне разбиения, меньшая плитка прикладывается левым верхним углом, то есть и при последующих разбиениях ребер не появится.

Теперь заметим, что при первом разбиении макроплитки в правом нижнем углу появляется ребро, но обе макроплитки разбиения примыкают к правому нижнему углу своим левым верхним углом. То есть, при дальнейших разбиениях новых ребер не появится. К левому нижнему углу после первого разбиения будут примыкать две плитки разбиения, правым верхним и правым нижним углами, то есть тут появится еще одно ребро при втором разбиении, и больше их не будет.

Итак, внутри каждого угла каждой макроплитки будет проведено, максимум, два ребра.

\medskip

Разбиение плиток, начиная с некоторого момента, перестает менять степень выбранной вершины.
\end{proof}

\medskip

{\bf Следствие.} 1. Каждая вершина заданной глубины $x$ выступает в качестве ядра подклейки для ограниченного количества вершин.

2. В каждую вершину входит ограниченное количество ребер различных уровней, включая ребра из подклеек. 

Первое утверждение следует из того, что каждая вершина бывает ядром подклейки только
один раз, когда ее глубина на $2$ меньше максимальной глубины.  Докажем второе. Для выбранной вершины из ее собственной
 макроплитки ребер идет конечное число. 
Ребра одного уровня проводятся в вершину при одновременном создании нескольких подклеек, где данная вершина лежит на 
ребре подклейки. Таких подклеек одновременно
создается конечное число, так как ядра всех этих подклеек находятся на расстоянии $2$ от нашей вершины.

 \medskip
 
  {\bf Нумерация на входящих в вершину ребрах.} Для каждой вершины пронумеруем входящие в нее ребра. Нумеруются именно способы входа в вершину, то есть в случае, если вершина лежит на некоторой стороне или внутреннем ребре макроплитки, то эта сторона или ребро даст два входящих в вершину ребра.
  
    Будем делать это сначала для ребер вершин глубины 
  $1$ и $2$. У таких вершин нет ребер в подклейки. Все плоские ребра можно пронумеровать по уровням. Сначала даем номера всем ребрам самого большого уровня, потом на один меньшего и так далее. При равных уровнях нумеруем по часовой стрелке.
  
   Допустим, все входящие ребра в вершины глубины $k$ пронумерованы.
  Пронумеруем все ребра из вершины $X$ глубины $k+1$. Все входящие ребра в ту же базовую плоскость можно опять пронумеровать по уровням, от большего к меньшему, а при равных уровнях -- по часовой стрелке. Рассмотрим все ребра, выходящие из $X$ в подклееные области. Пронумеруем сначала ядра этих областей, от меньшей глубины к большей, при равной глубине сначала боковые и краевые вершины, потом внутренние типа $C$, потом $A$ и $B$. Заметим, что $X$ не может лежать в двух подклейках ядра которых имеют одну и ту же глубину и  одновременно являться боковыми или краевыми, а также внутренними одного типа. Действительно, в этом случае две подклейки возникают одновременно с разными ядрами глубины $n-2$ при максимальной глубине $n$. В этом случае расстояние между ядрами в момент подклейки будет не более $4$. В этом случае ядра не могут одновременно быть внутренними вершинами одного типа или одновременно быть боковыми или краевыми.
  
  Теперь покажем, как пронумеровать ребра для подклеек с одним и тем же ядром.
  Рассмотрим некоторое ребро $r$ из $X$ в подклееную область $T$. Эта область (макроплитка) $T$ 
  подклеивается к остальному комплексу
  по двум своим сторонам, причем $X$ лежит на одной из них. Отметим, что $X$ не может быть ядром данной подклейки, так как ядро по отношению к подклейке является верхней левой вершиной, а из верхней левой вершины макроплитки не выходит внутренних ребер.
  
   Две стороны, по которым подклеивается $T$,  являются входящими ребрами в ядро данной подклейки $T$.
  У ядра подклейки глубина меньше, и нумерация на ребрах уже задана. 
  Поставим в соответствие ребру $r$ номер ребра, соответствующего стороне, на которой $X$ {\it не} лежит. 
 Совпадение номеров у ребер $r_1$ и $r_2$ означает, что они идут в одну и ту же подклейку. Фактически мы пронумеровали все подклейки,в которых участвует вершина $X$.

  Рассмотрим теперь все ребра в подклейки из $X$, имеющие некоторый одинаковый
  уровень $s$. Пронумеруем все ребра, соответствующие одинаковому номеру подклейки ( входящего в соответствующее ядро ребра, на котором $X$ не лежит).  Заметим, что выходящие ребра из вершины на границе макроплитки, если они имеют один уровень, то тогда у них разный тип ребра. То есть у таких ребер разный тип ребра и их можно пронумеровать
  от меньшего к большему. Далее так нумеруем все ребра для каждого номера подклейки от меньшего к большему.
  
Проделав эту операцию по очереди для уровней ребра $s$ от большего к меньшему можно пронумеровать все ребра в выделенную вершину любой заданной глубины.
  
  \medskip
 
 % Ребра максимального уровня, входящие в вершину, будем называть {\it  главными}.  

%Плоские части комплекса представляют собой макроплитки некоторого уровня. Каждая такая макроплитка либо является начальной, базовой плоскостью, либо в какой-то момент была подклеена, а потом несколько раз разбита. Макроплитки состоят из плиток (простых циклов длины $4$).

\medskip

\begin{lemma}[О боковой вершине] \label{side_node}

1) Каждая боковая вершина лежит на середине стороны в какой-либо макроплитке или в двух макроплитках одного уровня, лежащих в одной плоскости;

2) Если боковая вершина не находится на границе исходного комплекса первого уровня, то она лежит на одном из восьми внутренних ребер (рисунок~\ref{fig:inneredges}) в некоторой макроплитке, либо лежит на границе подклееной макроплитки.

\end{lemma}

\begin{proof}
В комплексе первого уровня свойство выполняется. На втором уровне также все верно для всех созданных боковых вершин.

Пусть для уровня $k$ комплекса все эти свойства выполнены. Для уже существующих боковых вершин оба свойства будут сохраняться при дальнейших разбиениях и подклейках. Если новая боковая вершина возникла при подклейке, то это она использовалась в качестве $T_2$ или $T_3$ из определения подклейки, и в качестве макроплитки можно взять как раз подклееную в тот момент макроплитку.

 Созданная при разбиении вершина, очевидно, лежит на середине стороны только что разбитой плитки. Если это подклееное ребро или граница начального комплекса первого уровня, то такая плитка одна (и тогда это граница подклееной макроплитки, либо всего комплекса), в остальных случаях таких макроплиток две.

Для уровней $1$ и $2$ комплекса заметим, что все стороны плиток, не лежащие на границе, являются частью какого-то большего ребра, классифицируемого, как одного из восьми типов из рисунка~\ref{fig:inneredges}. Заметим, что если боковая вершина лежит на внутреннем ребре определенного типа, то при дальнейших подклейках и разбиениях этот тип не меняется. При этом создаваемые при разбиении внутренние ребра, очевидно, принадлежат одному из этих восьми типов.

Если ребро, при разбиении которого образовалась вершина не является границей для всего комплекса или подклееной макроплитки, то  оно принадлежит одному из восьми типов из рисунка~\ref{fig:inneredges} для какой-то макроплитки. Значит и вершина тоже принадлежит такому ребру.

\end{proof}

\medskip

%{\bf Примечание.} Для каждой боковой вершины $X$ можно выделить пару вершин $(Y,Z)$, являющихся концами ребра, на котором лежит $X$. Эту пару вершин будем называть {\it ``начальниками''} для $X$.

\medskip

\subsection{Пути на комплексе}

В этом параграфе мы установим несколько свойств путей, проходящих по комплексу. В леммах ниже мы будем рассматривать комплексы произвольного размера, то есть мы дополнительно полагаем, что они верны для комплекса произвольного размера, естественно, на котором могут существовать указанные в леммах конструкции. Если речь идет о существовании некоторых констант, достаточно длинных путей и т.п. мы будем считать что такие константы и прочие выборы производятся независимо от размеров комплекса.

\medskip

Путь, идущий по макроплитке, может быть локально преобразован:

\begin{lemma}[О переброске пути]  \label{path_flip}

Пусть $XYZT$ -- некоторая макроплитка. Рассмотрим путь $XYZ$ (состоящий из двух соседних макроребер). Тогда, если разрешается менять подпуть из двух соседних ребер любой плитки на подпуть из двух других ребер (с общими началом и концом у этих подпутей), то путь $XYZ$ может быть преобразован в $XTZ$ -- путь по другим двум макроребрам.

\end{lemma}

\begin{figure}[hbtp]
\centering
\includegraphics[width=0.7\textwidth]{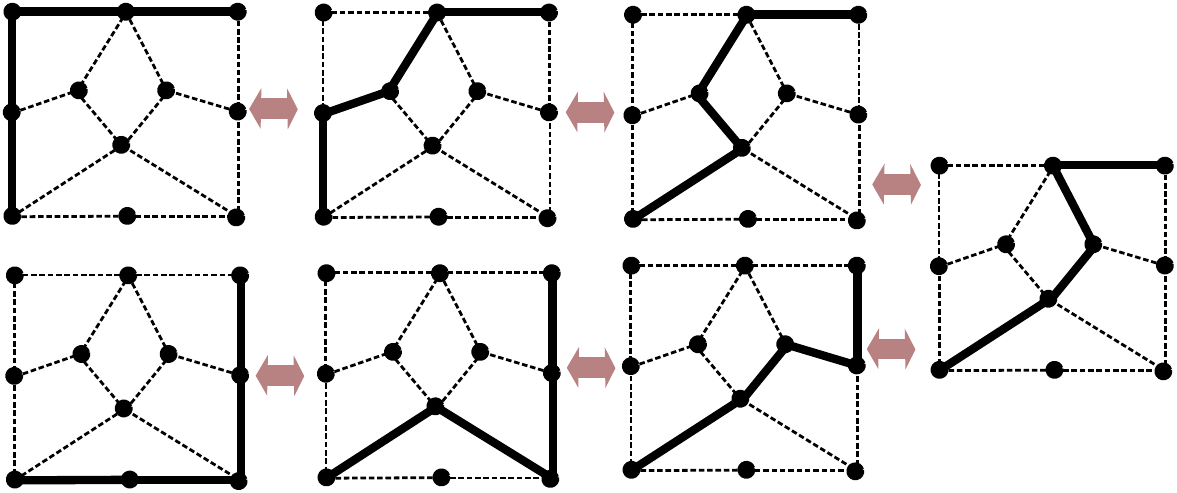}
\caption{Переброска пути на другую сторону макроплитки.}
\label{fig:transformpath}
\end{figure}

\begin{proof} Будем доказывать по индукции. Для случая, когда макроплитка -- это плитка, преобразование можно сделать сразу. Шаг индукции можно совершить, выполняя
локальные преобразования по правилам, показанным на рисунке~\ref{fig:transformpath}.

\end{proof}

\medskip

{\bf Примечание.} Применяя аналогичные рассуждения, можно доказать то, что также перебросить можно и путь, начало и конец которого лежат в серединах противоположных ребер макроплитки. ,

\medskip

\begin{definition}
Замены подпутей из пары соседних ребер плитки на подпуть из другой пары будем называть {\it локальными заменами}. А возможность с помощью таких замен достичь некоторого состояния, будем называть {\it локальным преобразованием}.
\end{definition}

\medskip

\begin{definition}
 Будем считать форму пути {\it нулевой}, если она содержит подпуть длины $2$ по ребру некоторой минимальной плитки (туда-обратно).
\end{definition}

\medskip

\begin{lemma}[О выносе пути на границу]  \label{to_border}

Пусть начало и конец пути $P$, проходящего по макроплитке $T$, лежат на границе $T$. Тогда можно реализовать одну из двух возможностей:

1) $P$ может быть локально преобразован в нулевую форму.

2) $P$ может быть локально преобразован в форму $P'$ так, что $P'$ полностью лежит на границе $T$.

Кроме того, любой кратчайший путь, соединяющий противоположные углы или середины противоположных сторон макроплитки имеет длину $2^n$, и может быть локально преобразован в любой из двух путей по границе (полупериметр), с теми же концами.

\end{lemma}

\begin{proof} Будем доказывать индукцией по уровню макроплитки.
Любой путь по минимальной плитке и так лежит на ее границе. Пусть макроплитка уровня $n$. Рассмотрим ее разбиение на $6$ (макро)плиток уровня $n-1$. Путь $P$ разбивается на несколько участков, каждый из которых имеет начало и конец на границе какой-то из этих подплиток. По предположению индукции, эти участки можно локально преобразовать в форму $P_1$ так, чтобы весь путь $P_1$ проходил только по границам подплиток.

\medskip

{\bf Замечание.} Символы $L$, $U$, $R$, $D$, $A$, $B$, $C$ обычно обозначают типы вершин. Соответственно, внутри макроплитки крупного уровня может быть несколько вершин, имеющих один и тот же тип. Но в данном предложении макроплитка фиксирована, так что логично для обозначения принадлежащих ей вершин использовать те же буквы.

\medskip

Обозначим буквами $L$, $U$, $R$, $D$ середины соответственно левой, верхней, правой, нижней сторон $T$, и внутренние вершины как $A$, $B$, $C$.

%Допустим, путь $P_1$ содержит дважды некоторую вершину  нашей макроплитки $T$. Тогда он должен полностью обходить некоторую подмакроплитку уровня $n-1$, примыкающую к этой вершине. В этом случае можно применить лемму о переброске пути и получить форму, где выполнено 1).

Пусть $P_1$ содержит $A$ один раз. Войти и выйти $P_1$ может только по трем ребрам (направленным в сторону $U$, $L$ и $C$). Если вход и выход -- по одному ребру, то условие 1) выполнено.

Пусть во входе и выходе задействованы ребра в сторону $U$ и $L$. Если, например, выйти в сторону $L$, то сойти с ребра мы не можем (так как $P_1$ проходит только по границам макроплиток $n-1$ уровня) и вернуться назад тоже (так как тогда будет выполнено 1). Заметим, что тогда $P_1$ содержит подпуть $LAU$ или $UAL$, в обоих случаях этот подпуть, согласно лемме о переброске пути, можно локально преобразовать чтобы он проходил по границе макроплитки $T$ и, таким образом, не содержал $A$.

Если во входе и выходе задействованы ребра в сторону $L$ и $C$, то $P_1$ содержит подпуть $LAC$ или $CAL$ и он может быть преобразован в форму соответственно $LXC$ или $CXL$, где $X$ -- левый нижний угол, и в этом случае опять форма не будет содержать $A$.

Если во входе и выходе задействованы ребра в сторону $U$ и $C$, то $P_1$ содержит подпуть $UAC$ или $CAU$ и он может быть преобразован в форму соответственно $UBC$ или $CBU$, и тоже форма не будет содержать $A$.

\medskip

Проделав аналогичные рассуждения для симметричного случая вершины $B$, мы можем заключить, что путь можно преобразовать в форму $P_2$, такую что она либо вообще не содержит $A$ и $B$, либо содержит кусок $UBC$ или $CBU$ (который может быть преобразован в $UAC$ или $CAU$ соответственно).

Если $P_2$ не содержит $A$ и $B$ и  содержит $C$, то войти и выйти из $C$ можно только по ребрам, ведущим в левый нижний и правый нижний углы. В этом случае можно применить лемму о переброске пути для участка, содержащего эти углы и $C$ и получить форму, не содержащую вершин $A$, $B$, $C$.

Пусть $P_2$ содержит кусок $UBC$ или $CBU$. Посмотрим, как $P_2$ проходит через $C$: один вход (или выход) в сторону $B$, другой не может идти в сторону $A$, так как $P_2$ не содержит $A$, значит он идет в сторону нижнего левого или нижнего правого угла. Если это нижний правый угол $Y$ то $P_2$ содержит кусок $UBCY$ или $YCBU$. Применяя лемму о переброске путей, $UBCY$ переводится в $UBRY$ и далее в $UZRY$, где $Z$ правый верхний угол. Аналогично, $YCBU$ переводится в $YRZU$.
Если это нижний левый угол $X$, то $P_2$ содержит кусок $UBCX$ или $XCBU$. Аналогично применяя лемму~\ref{path_flip} о переброске пути, переводим $UBCX$ в $UACX$ и далее в $UALX$ и потом в $UWLX$, где $W$ -- левый верхний угол. Аналогично $XCBU\Rightarrow XCAU \Rightarrow XLAU   \Rightarrow XLWU$.

Таким образом, мы получим форму не содержащую вершин $A$, $B$, $C$. В этом случае, если нельзя выполнить условие 1), то весь путь проходит по границе $T$, и выполнено условие~2).

\medskip

Последнее утверждение этой леммы является следствием из первого, с учетом леммы о переброске пути.
\end{proof}

\begin{lemma}[О ``шевелении'' пути] \label{path_changing}

 Пусть $P$ -- путь, состоящий из макроребра $XY$ некоторой макроплитки $T$ уровня не менее $4$.
Обозначим середину $XY$ как $H$, середины $XH$ и $HY$ как $G_1$ и $G_2$.
 В соответствии с определением подклеек, существует подклееная макроплитка, три угла которой лежат на макроребре: это $G_1$, $H$ и $G_2$.
Тогда $P$ может быть локально преобразован в форму $P'$, состоящую из трех частей:
первая -- $XG_1$, третья -- $G_2Y$, а вторая представляет собой путь по двум соседним ребрам подклееной макроплитки, два других ребра которой являются частью $P$.
\end{lemma}

\begin{proof} Достаточно проследить разбиение $T$ (на три уровня в глубину).
Обозначим середины $G_1H$ и $HG_2$ как $G_3$ и $G_4$. Так как $G_3$ и $G_4$ будут иметь уровень на $2$ больше, чем $H$, значит, была проведена подклейка для пути $G_1G_3HG_4G_2$. Пусть $F$ -- созданная при этой подклейке угловая вершина. По лемме о переброске пути, $G_1HG_2$ можно локально преобразовать в $G_1FG_2$.

\end{proof}

\begin{lemma}[О выделении локального участка] \label{longpath}

Пусть $n$ -- уровень макроплитки и путь $P$ лежит внутри нее.

1. Пусть оба края $P$ лежат на границе $T$. Тогда если путь $P$ имеет длину не менее $5\times 2^{n-2}$, он может быть локально преобразован в нулевую форму $P'$.

2. Пусть один край $P$ лежит в углу или на середине стороны $T$, а второй край - внутри $T$ (часть пути может проходить по границе $T$). Тогда если путь $P$ имеет длину $2^{n}$, он может быть локально преобразован в нулевую форму $P'$.

\end{lemma}

\begin{proof}
1. Согласно лемме о выносе на границу, можно считать, что весь путь $P$ проходит по границе $T$. Пусть при этом не образовалось нулевой формы, то есть путь покрывает какую-то часть границы $T$, совершая движение либо по часовой стрелке либо против.  Отметим на границе $T$ восемь точек -- углы и середины сторон. Заметим, что расстояние между соседними отмеченными точками (углом и серединой стороны) равно $2^{n-2}$. Наш путь включает в себя не менее $5$ из отмеченных точек. Тогда найдутся две отмеченные точки, диаметрально противоположные относительно $T$, и такие, что путь $P$ полностью покрывает половину границы между ними. Заметим при этом, что хотя бы одна из этих отмеченных точек не является концевой для $P$. Пользуясь леммой о переброске (с примечанием) можно перевести подпуть из половины границы на другую половину.
Ясно, что получившаяся форма будет нулевой.

\medskip

2.Будем доказывать индукцией по уровню макроплитки. Для уровней $1$ и $2$ утверждение легко проверяется. Допустим, уровень макроплитки равен $n$, а для меньших $n$ это верно. Будем размечать пути скобками: открывающей скобкой отметим вершину, где путь заходит внутрь некоторой макроплитки, а закрывающей - где выходит. Из всех локально эквивалентных пути $P$ путей выберем такой путь $Q$, в котором меньше всего скобок. Заметим, что подпуть $Q$ , начинающийся с открывающей скобки, заканчивающийся закрывающей и не содержащий скобок внутри, попадает под условия леммы о выносе пути на границу, и его можно преобразовать в форму, где скобок будет меньше. Также заметим, что один край $Q$ на границе, то есть первая скобка (если она есть) обязана быть открывающей. Тогда в $Q$ не может быть закрывающих скобок, только открывающие.

\smallskip

Если $Q$ не содержит скобок, тогда он проходит по границе макроплитки, и тогда оба конца пути лежат на границе, что не так по условию.

%Путь $Q$ входит внутрь некоторой макроплитки $T'$ ($T'$ может совпадать с $T$) далее до второй скобки идет по внутренним ребрам, разделяющим $T'$ на шесть ее дочерних макроплиток. Если второй скобки вообще нет, то длина $Q$ не может превышать $2^n-1$. Пусть она есть. Тогда на второй скобке путь заходит внутрь некоторой макроплитки $T''$. Заметим, что либо $T''$ является одной из шести дочерних макроплиток $T'$, либо лежит внутри одной из них. Обозначим такую дочернюю макроплитку $T'$ как $W$. Тогда путь $Q$ можно разбить на участки $Q_1$ и $Q_2$, где $Q_2$ начинается с первой вершины, лежащей на границе $W$. (эта вершина может быть помечена второй скобкой или встретиться в пути ранее).  Заметим, что $Q_1$ идет по внутренним ребрам $T'$, причем его длина не превышает $2\cdot 2^{n-2}$, так как в макроплитке второго уровня расстояние от края до любой макроплитки не более $2$ и длина одного внутреннего макроребра в $T'$ не более $2^{n-2}$.  $Q_2$ удовлетворяет условиям предположения индукции, и его длина не превышает $2^{n-1}$. В случае 2 переход завершен.

Обозначим как $H$ одну из шести подплиток $T$, в которой лежит конец пути $Q$. Заметим, что если путь заходит в какую либо подплитку, то выйти из нее он уже не может, так как в этом случае в пути будет закрывающая скобка. 

%Пусть начало $Q$ тоже лежит в $H$. Так как путь не может зайти в подплитку, отличную от $H$, то он весь лежит в $H$ и можно применить предположение индукции.

% Пусть начало не лежит в $H$. 
 
 Обозначим как $X$ первую вершину на пути, что вся оставшаяся часть лежит в $H$. Ясно, что $X$ -- один из углов одной из шести подплиток $T$. Разобьем $Q$ на два подпути $Q_1$ и $Q_2$, до и после $X$. Заметим, что
$Q_1$ проходит лишь по границам и по внутренним ребрам $T$.  Таким образом, можно считать, что $Q_1$ проходит по внутренним ребрам и границам макроплитки второго уровня.

Теперь проверим, что на макроплитке уровня $2$ любой путь $Q_1$, начинающийся на границе, ведущий в подплитку $H$ и не содержащий закрывающих скобок, приводится либо к нулевой форме, либо к виду $S_1S_2$, где $S_1$ имеет длину не более $2$, а $S_2$ лежит в $H$. В силу вышесказанного, достаточно рассмотреть случаи, когда $Q_1$ имеет длину $3$ и оканчивается во внутренней вершине, либо имеет длину $4$ и лежит на границе макроплитки. Требуемое утверждение проверяется перебором путей, подчиняющихся указанным ограничениям, и подплиток $H$, в которых эти пути оканчиваются. Детали перебора предоставляются читателю.»

%Заметим, что из любой середины стороны или угла макроплитки второго уровня можно добраться до какого-нибудь угла любой выбранной подплитки пройдя не более, чем по двум сторонам подплиток $T$. 

Итак, можно считать, что начальная часть $S_1$ пути  не длиннее $2\times 2^{n-2}$. Тогда оставшаяся часть пути лежит в $H$, и в ней можно выбрать кусок длины $2^{n-1}$. К нему можно применить предположение индукции.

%Если $Q_1$ сам по себе не длиннее $2\times 2^{n-2}$, то предположение индукции можно применить к $Q_2$.

%Пусть $Q_1$ длиннее $2\times 2^{n-2}$. Итак, путь $Q_1$ проходит по сторонам и внутренним ребрам макроплитки второго уровня, начинается в углу или середине стороны, заканчивается в углу некоторой подплитки $H$. 

%Покажем, что такой путь может быть локально преобразован к форме $S_1S_2$, где $S_1$ -- путь не более чем из двух сторон подплиток $T$, а $S_2$  проходит по границе $H$.

\smallskip

\end{proof}

\begin{lemma}[О непродолжаемом пути] \label{bad_path}

Рассмотрим некоторую макроплитку $T$, являющуюся подплиткой более крупной макроплитки.
Пусть путь $P$ лежит в $T$, причем начало и конец $P$ лежат в серединах противоположных сторон $T$. Тогда любые плоские пути вида $WP$ или $PW$, где длина $W$ более длины $P$, могут быть приведены к нулевой форме.
\end{lemma}

\begin{proof} Случаи $WP$ и $PW$ аналогичны, рассмотрим $PW$. Пусть $a$ -- это та сторона $T$, в середине $X$ которой заканчивается $P$ и начинается $W$. Пусть путь, выйдя из $X$ идет внутрь макроплитки $T$, либо по ее стороне $a$. Обозначим как $r$, его первое ребро. Если $r$ уходит внутрь макроплитки $T$, то к пути $Pr$ применима вторая часть леммы о выделении локального участка. Пусть $r$ лежит на стороне. Тогда один из двух путей $Q_1$ или $Q_2$ c теми же концами, что и $P$ (и локально эквивалентных ему по лемме о выносе пути на границу) тоже будет содержать $r$. В этом случае $Q_1r$ или $Q_2r$ содержит ребро, которое проходится туда-обратно.

%Заметим, что существуют лежащие на границе $T$ пути $Q_1$ и $Q_2$, с теми же концами, что и $P$, обходящие $T$ по границе, с разных сторон. Заметим, что по леммами о выносе путей на границу и о переброске пути, $P$ локально преобразуется в $Q_1$, и $Q_1$ в $Q_2$. Рассмотрим множество $M$ всех последних ребер (заканчивающиеся в $X$) всех путей-промежуточных этапов преобразования $Q_1$ в $Q_2$. $М$ образует <<пучок>> -- несколько ребер с общей вершиной $X$, причем $Q_1$ и $Q_2$ соответствуют крайним и противоположным ребрам $M$, идущим по границе $T$. Ясно, что при одном локальном преобразовании последнее ребро либо не меняется, либо меняется на соседнее ребро в $М$. Поскольку $r$ содержится в $M$, существует некоторый эквивалентный $Q_1$, $Q_2$ и $P$ путь $P'$ с последним ребром $r$.  Но тогда путь $P'r$ содержит ребро $r$, которое проходится туда-обратно.

\smallskip

Значит, из $X$ путь может выйти только в макроплитку $T'$, соседствующую с $T$ по стороне $a$. Допустим, существует вершина $Y$ из пути $W$, лежащая на границе $T'$, отличная от $X$. В этом случае, кусок пути $W$ от $X$ до $Y$ может быть локально преобразован так, что он будет проходить по границе $T'$. Но тогда выход этого куска из $X$ будет по стороне $a$, и в этом случае опять можно локально преобразовать $P$ чтобы $PW$ содержал туда и обратно по одному ребру.  Если такой вершины $Y$ не найдется, то все остальные вершины $W$ лежат внутри $T'$, то есть можно воспользоваться леммой о выделении локального участка.

\end{proof}

\smallskip

\begin{definition}
Пусть $r$ -- выходящее ребро из вершины $X$. Пусть $X$ -- боковая вершина и разделяет некоторые макроплитки $U$ и $W$. $r$ называется {\it главным ребром}, если оно лежит на границе между $U$ и $W$. Если же $X$ -- внутренняя вершина в макроплитке $U$, тогда $r$ {\it главное ребро}, если $r$ -- лежит на одном из внутренних ребер макроплитки $U$.
\end{definition}

\begin{definition}
Рассмотрим некоторый плоский путь. Он проходит через внутренние и боковые вершины. Рассмотрим последовательность всех вершин, которые посещает путь. Если на пути есть некоторая боковая вершина, вход и выход в которую прошел по главным ребрам, выбросим ее из этой последовательности. Последовательность оставшихся вершин назовем {\it паттерном пути}. То есть в паттерн пути входят все внутренние вершины, и все боковые вершины, где происходит вход(выход) в (из) макроплитку.
Фактически паттерн пути -- это его карта, показывающая маршрут на ориентировочных точках.
\end{definition}

\begin{definition}
Пусть $P$ -- некоторый паттерн. В рамках данного определения путь $W$ будем называть {\it достаточно большим} путем, содержащим $P$, если выполнены следующие условия:

1) Путь $W$ представляется в виде $W_1W_2W_3$, где паттерн $W_2$ есть $P$;

2) Длина как $W_1$, так и $W_3$ не менее удвоенной длины $W_2$.

{\it Мертвым} будем называть такой паттерн $P$, что любой достаточно большой путь $W$, cобственным подпутем которого является путь с паттерном $P$, может быть локально преобразован к нулевой форме.
\end{definition}

\begin{lemma}[О мертвых паттернах] \label{DeadPaterns}

Рассмотрим некоторую макроплитку $T$ и обозначим в ней внутренние вершины $A$, $B$, $C$ и боковые $U$, $R$, $D$, $L$ (аналогично обозначениям при разбиениях).
Тогда паттерны $AUB$, $ACB$, $CXD$ (где $X$ -- нижняя левая, либо нижняя правая вершина) являются мертвыми.
\end{lemma}

\begin{proof} Пути с паттернами $AUB$ и $ACB$ локально преобразуются друг в друга, так что достаточно рассмотреть один из этих случаев. Пусть $n$ -- уровень $T$. Пусть $W_1ACBW_2$ -- достаточно большой путь, содержащий подпуть с паттерном $ACB$. Заметим, что длина пути с паттерном $ACB$ равна $2^{n-1}$. Допустим, $W_1$ не содержит вершин на границе $T$.

Если на $W_1$ есть вершины, лежащие на внутренних ребрах $T$, то участок пути от первой такой вершины $K$ до $A$ может быть локально преобразован так, чтобы он проходил только по внутренним ребрам $T$. Заметим, что таким внутренним ребром может быть только ребро от $A$ до середины левой стороны, в других случаях путь $W_1ACB$ либо будет содержать участок, которое он проходит туда, а потом обратно, либо будет содержать подпуть по границе средней подплитки $T$ длиной более $2^{n-1}$. То есть может быть преобразован к нулевой форме.

Тогда можно считать, что $W_1$ состоит из двух участков $W_11$ и $W_12$, причем $W_11$ проходит внутри левой верхней или левой нижней подплитки $T$, а $W_12$ -- по ребру от середины левой стороны $T$ до $A$. Длина $W_1$ не менее $2^{n-1}$. Применяя лемму о выделении локального участка, получаем, что путь $W_1$, заканчивающийся в $A$ приводится к нулевой форме.

Итак, пусть теперь $W_1$ и, аналогично, $W_2$ содержат вершины на границе $T$. Рассмотрим разбиение макроплитки $T$ на $6$ дочерних подплиток (в соответствии с операцией разбиения) и будем локально преобразовывать $W$ так, чтобы он проходил по границам этих дочерних подплиток. Тогда из вершины $B$ путь
должен идти по ребру к $R$, а в $A$ путь должен входить по ребру из $L$, в остальных случаях образуется обход подплитки по трем ее сторонам (и этот кусок преобразуется к нулевой форме). Таким образом, в $W$ можно выделить подпуть, проходящий по макроплитке $T$, начало которого в $L$ и конец в $R$. Можно воспользоваться леммой о непродолжаемом пути и получить требуемое.

Для паттерна $CXD$, пусть $T'$ -- макроплитка, соседняя с $T$ по нижней стороне. Заметим, что выход из вершины $D$ обязательно должен быть внутрь макроплитки $T'$, иначе можно применить лемму о выделении локального участка для нижней дочерней подплитки $T$. Заметим, что длина пути с паттерном $CXD$ равна $2^{n-1}$, где $n$ -- уровень $T$ (и $T'$).  Дальнейшая часть пути не менее, чем вдвое длиннее, то есть, не менее $2^{n}$. Если далее на пути не встречается вершин на границе $T'$, то можно применить лемму о выделении локального участка. Пусть на пути встретится вершина $Z$ на границе $T'$.  Преобразуем  кусок $DZ$ чтобы он проходил по периметру $T'$.

Но тогда из $D$ путь $W$ продолжится по граничной между $T$ и $T'$ стороне.
То есть теперь опять можно применить лемму о переброске пути для нижней дочерней подплитки $T$ и получить требуемую нулевую форму.

\end{proof}

\medskip

{\bf Примечание.} Ясно, что и паттерны $BUA$, $BCA$, $DXC$ (где $X$ -- нижняя левая, либо нижняя правая вершина) тоже мертвые, доказательство полностью аналогично.

\medskip
\begin{lemma}[О мертвых путях в нижней подплитке] \label{DeadPaths}

Рассмотрим некоторую макроплитку $T$ уровня $n$. Пусть путь $XYZ$ лежит в $T$, причем $Y$ -- левый нижний угол $T$, $XY$ лежит на внутреннем ребре, идущем из левого нижнего угла $T$ во внутреннюю вершину $C$, а $YZ$ лежит на нижней стороне $T$. Тогда для любых плоских путей $W_1$, $W_2$, длины которых более $2^{n+1}$, путь $W_1XYZW_2$ можно локально преобразовать к нулевой форме.

\end{lemma}

\begin{proof} Обозначим середину нижней стороны $T$ как $D$. Пусть $T'$ нижняя подплитка $T$.  Рассмотрим некоторый плоский путь $W_1XYZW_2$, c длинами $W_1, W_2$ большими $2^{n+1}$. Учитывая длину путей  $W_1$ и $W_2$, они имеют вершины на границе макроплиток уровня $n$. Тогда их можно преобразовать так, чтобы они проходили по границам макроплиток уровня $T'$ или выше (не заходя в более мелкие плитки). Ясно, что $W_1$ должен проходить через $C$, а $W_2$ через $D$, иначе возникает очевидный кусок с нулевой формой. Таким образом, после преобразования, путь будет содержать подпуть с паттерном $CXD$,  длина которого $2^{n-1}$ и который является мертвым, по лемме о мертвых паттернах. Следовательно, для указанных $W_1$ и $W_2$ наш путь можно преобразовать к нулевой форме.

\end{proof}

{\bf Замечание}. Доказательство аналогично переносится на случай, когда $Y$ -- правый нижний угол и $XY$ идет по внутреннему ребру от $C$.

\medskip

\begin{definition}
 {\it Некорректным участком} пути будем называть такой подпуть $XYZ$, что вершина $Y$ лежит на границе некоторой макроплитки $T$, а вершины $X$ и $Z$ лежат внутри $T$ (не на границе).
\end{definition}

\begin{lemma}[О некорректных участках] \label{uncorrect_sectors}

Пусть есть некорректный участок $XYZ$ в макроплитке $T$ уровня $n$, причем $T$ -- минимальная макроплитка, содержащая $XYZ$ в качестве некорректного участка. Тогда для любых плоских путей $W_1, W_2$, длины более $2^{n+2}$, путь $W_1XYZW_2$, может быть локально преобразован к нулевой форме.
\end{lemma}

\begin{proof}
Пусть $U$ - родительская макроплитка для $T$.
Пусть $W_1$ не содержит вершины, лежащей на границе $U$.
В силу минимальности выбора $T$, вершина $Y$ лежит в углу, либо в середине стороны $T$.
Если на $W_1$ нет вершины, лежащей на границе $T$, тогда к пути $W_1XY$ можно применить лемму о выделении локального участка. Значит, на пути $W_1$ есть вершины, лежащие на внутренних ребрах $U$ (в частности, на границе $T$). Пусть $K$ первая вершина в $W_1$, лежащая на внутреннем ребре $U$.
Тогда по лемме о выносе пути на границу, участок пути $W_1XY$ от $K$ до $Y$ может быть локально преобразован так, что он будет проходить только по внутренним ребрам $U$. Заметим, что тогда участок пути $W_1$ от начала до $K$ лежит в некоторой подплитке $U$, обозначим ее $T'$.

Допустим сначала, что $K$ и $Y$ лежат на одном внутреннем ребре (разделяющему $T$ и $T'$). Тогда после преобразования участка от $K$ до $Y$, весь путь от начала $W_1$ до $Y$ лежит в макроплитке $T'$ и заканчивается в середине стороны или в углу. То есть к нему применима лемма о выделении локального участка.

Пусть теперь $K$ и $Y$ лежат на разных внутренних ребрах. Тогда участок от $K$ до $Y$, проходящий по внутренним ребрам $U$ обязан будет пройти через один из концов внутреннего ребра, на котором лежит $K$. Пусть это вершина $F$ (учитывая, что $F$ не на краю $U$, это может быть лишь одна из трех внутренних вершин). Тогда длина участка пути он начала до $F$ меньше $2^n$, иначе к нему применима лемма о выделении локального участка. Тогда длина участка от $F$ до $Y$ более $3\times 2^n$, то есть более $6$ внутренних ребер $U$. В этом случае, очевидно, этот участок может быть преобразован к нулевой форме.

\medskip

Итак, будем считать, что $W_1$ и, аналогично, $W_2$ содержат вершины, лежащие на границе $U$.  Тогда, с учетом леммы о выносе пути на границу, пути $W_1$ и $W_2$ внутри $U$ проходят по внутренним ребрам, не менее крупным, чем ребра $XY$ и $YZ$. Кроме того, если вдоль пути встречается вершина на более крупном ребре, то и далее путь проходит по столь же крупным ребрам.

Как уже было сказано выше, вершина $Y$ лежит в углу $T$, либо в середине стороны. Из иерархии разбиений следует, что из всех углов кроме нижнего левого и из середины нижней стороны выходит не более одного ребра внутрь $T$, так что $Y$ может лежать только на середине левой, правой или верхней стороны или в левом нижнем углу.

{\bf 1.} Пусть $Y$ -- левый нижний угол. Из этого угла внутрь выходит два ребра, оба к внутренним C-вершинам, одна в макроплитке $T$, другая в нижней дочерней подплитке $T$. То есть, путь $XYZ$ попадает под условие леммы о нижней подплитке и все доказано.

{\bf 2.} Пусть $Y$ --  середина левой стороны. Из этой вершины выходит три ребра и поэтому есть три варианта расположения пути $XYZ$. В двух случаях можно применить лемму о нижней подплитке. Оставшийся случай изображен на рисунке~\ref{fig:badwayleftcase}.

\begin{figure}[hbtp]
\centering
\includegraphics[width=0.5\textwidth]{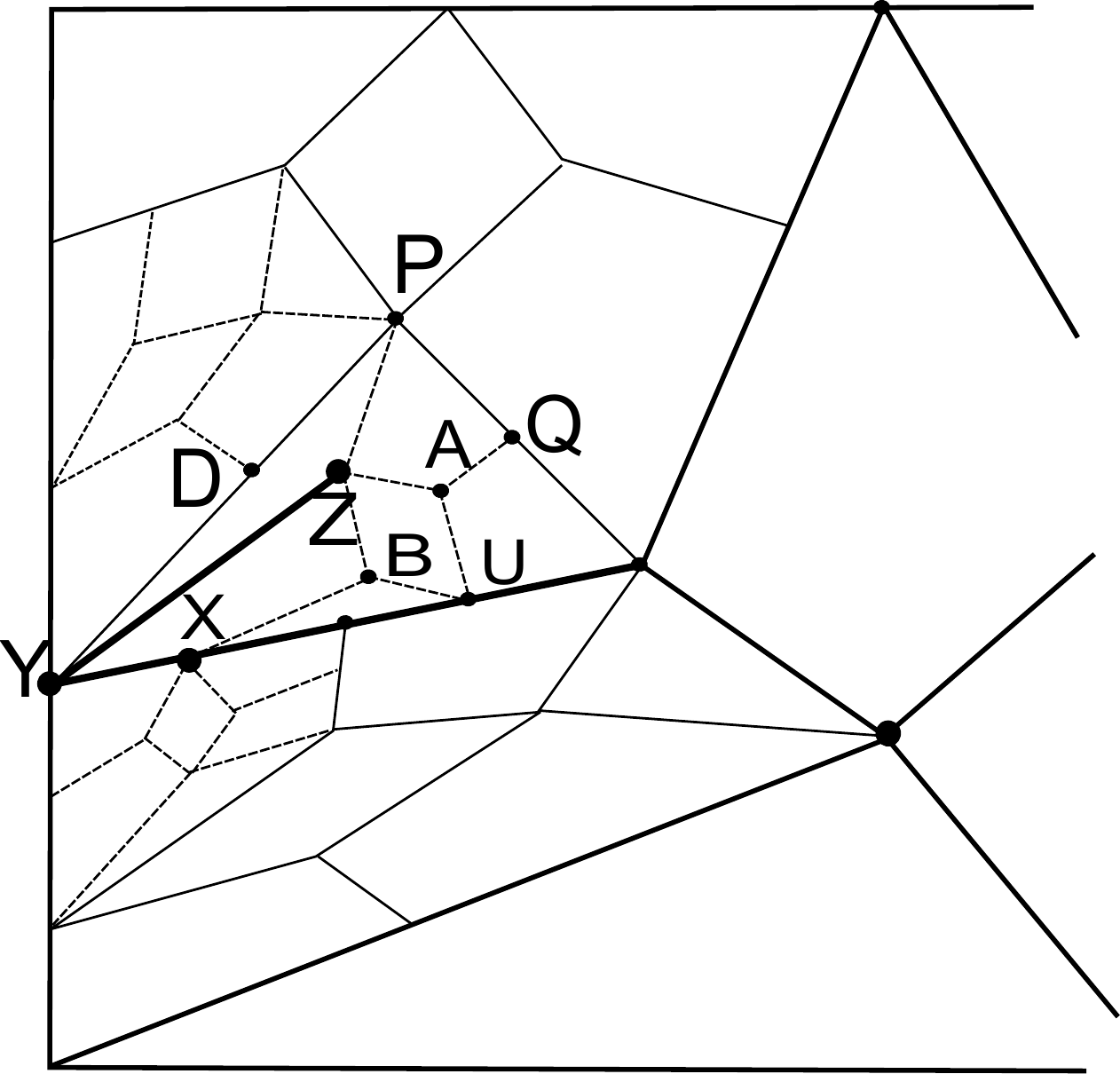}
\caption{Вершина $Y$ -- середина левой стороны.}
\label{fig:badwayleftcase}
\end{figure}

Из вершины $Z$ дальше путь может пойти:

 i) в вершину $B$. В этом случае кусок $BZYX$ по лемме о переброске пути преобразуется в $BXYX$, то есть в нулевую форму.

 ii) в вершину $P$. В этом случае кусок $YZP$ преобразуется в $YDP$ и для участка $XYD$ можно применить лемму о нижней подплитке.

 iii) в вершину $A$. Посмотрим, куда путь может пойти дальше. Если это вершина $Q$, то участок $YZAQ$ преобразуется в $YZPQ$ и далее в  $YDPQ$, после чего опять можно применить лемму о нижней подплитке. Во втором случае, если путь идет в вершину $U$, то путь $XYZAU$ сразу можно привести к нулевой форме:
$XYZAU\rightarrow XYZBU  \rightarrow XYXBU$.

\medskip

{\bf 3.} Пусть $Y$ --  середина правой стороны. Этот случай симметричен второму, только вместо левой верхней подплитки мы имеем дело с правой нижней. Все рассуждения полностью аналогичны второму случаю.

{\bf 4.}  Пусть $Y$ --  середина верхней стороны. Из этой вершины исходят четыре ребра
(рисунок~\ref{fig:badwayupcase}).

\begin{figure}[hbtp]
\centering
\includegraphics[width=0.6\textwidth]{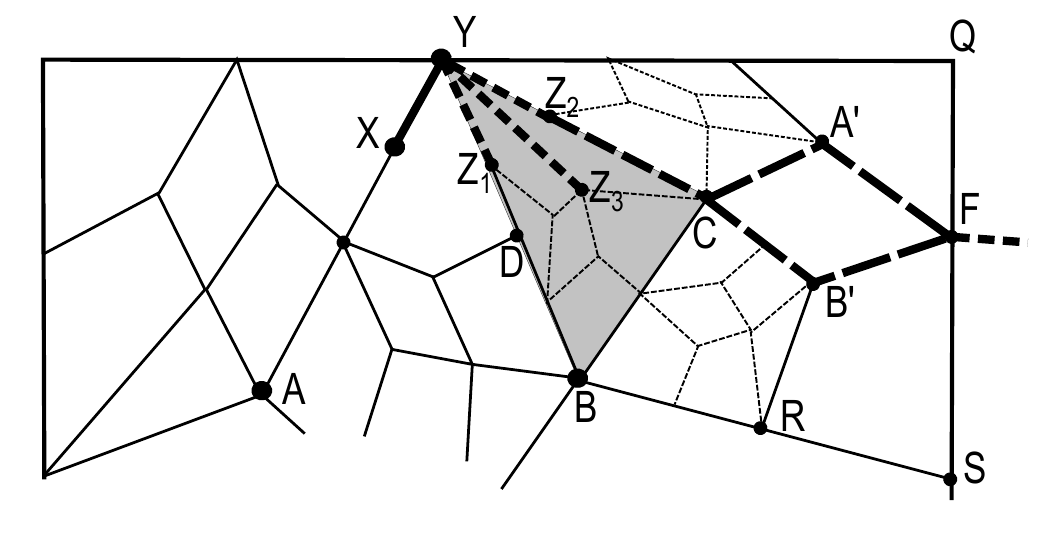}
\caption{Вершина $Y$ -- середина верхней стороны.}
\label{fig:badwayupcase}
\end{figure}

В случае, когда $X$ и $Z$ не лежат на внутреннем ребре типа $1$ (идущем из $Y$ к внутренней вершине $A$ макроплитки $T$), то получается ситуация, полностью аналогичная второму и третьему случаям выше, на этот раз вместо левой верхней подплитки мы имеем дело с правой верхней. Итак, пусть верхняя правая подплитка $T$ это $T'$ и, ради определенности, $X$ лежит на ребре, уходящем из $Y$ к внутренней вершине $A$ макроплитки $T$, а $Z$ -- на одном из трех ребер, попадающих в $T'$ (отмечены пунктиром на рисунке~\ref{fig:badwayupcase}).

Учитывая соглашение в начале доказательства, $W_1$ проходит по внутренним ребрам $T$ или более крупных макроплиток. Тогда либо $W_1XY$ содержит нулевую форму, либо $W_1$ проходит через внутреннюю вершину $A$ макроплитки $T$. Теперь рассмотрим путь $W_2$. Докажем, что этот путь может быть либо приведен к нулевой форме, либо локально преобразован в путь, начинающийся с макроребра $YB$.

В случае, если он начинается с $Z_1$, мы сразу получаем, что $W_2$ проходит через внутреннюю вершину $B$. В случае $Z_3$ путь идет по внутренним ребрам макроплитки, выделенной серым цветом на рисунке~\ref{fig:badwayupcase}. То есть он может быть локально преобразован, и будет проходить через $B$ либо через $C$. В случае $Z_2$ он также проходит через $C$.

В обоих этих случаях далее путь идет по одному из участков $YCA'F$ или $YCB'F$, указанных длинным пунктиром (либо $W_2$ пройдет через $B$, либо будет может быть приведен к нулевой форме).

Таким образом, $W_2$ теперь проходит через вершину $F$. 
Пусть $W_2$ покидает макроплитку $T$, уходя в соседнюю макроплитку по правой стороне (короткий пунктир). Тогда длина оставшейся части $W_2$ (после $F$) более чем $2^{n+2}-2^{n+1}=2^{n+1}$. При этом эта часть пути $W_2$ оказывается внутри макроплитки, такой же по уровню, что и $T'$ и  начинается в середине ее стороны. Тогда к этой части пути применима лемма о выделении локального участка и лемма о вынесении пути на границу, что сводит доказательство к рассматриваемым ниже случаям.

 Если из $F$ путь $W_2$ идет к вершине $Q$, то оба возможных куска $YCA'FQ$ и $YCB'FQ$ по лемме о выносе пути на границу приводятся к виду $YQFQ$, то есть к нулевой форме. Если путь из $F$ идет в $S$, то кусок $YCA'FS$ (или $YCBFS$) приводится к $YDBRS$ по лемме о переброске пути, применяемой последовательно для кусков этого пути.

Итак, $W_1$ проходит через внутреннюю вершину $A$ макроплитки $T$, а $W_2$ -- через внутреннюю вершину $B$.
Теперь получившийся путь имеет подпуть с мертвым паттерном $AUB$. Заметим, теперь, что оставшиеся участки путей $W_1$ и $W_2$ превосходят требуемые в определении мертвого паттерна. То есть мы можем применить лемму о мертвых паттернах и привести путь $W_1XYZW_2$ к нулевой форме.

\end{proof}

{\bf Примечание.} Рассуждения не сильно меняются, если рассмотренная макроплитка $T$ имеет уровень $3$ или меньше. В этом случае вершина $Z_1$ совпадает с $D$, а $Z_3$ с $C$.

\medskip

\begin{lemma}[О корректности путей] \label{correct_paths}

Пусть путь $P$ представляет собой проход по двум соседним сторонам некоторой макроплитки $T$. Тогда любые локальные преобразования не могут привести $P$ к нулевой форме, а также к форме, содержащей некорректный участок.

\end{lemma}

\begin{proof} Сначала покажем, что $P$ не приводится к нулевой форме. Для двух точек комплекса определено расстояние как наименьшая из длин соединяющих их путей. Докажем, что путь, проходящий по двум соседним сторонам макроплитки является кратчайшим для вершин в его концах. Действительно, для плитки уровня $n$ кратчайший путь можно вынести на периметр, применяя соответствующую лемму. Путь по периметру соединяющий противоположные точки имеет длину $2^n$.

Локальные преобразования не меняют длину пути. В случае, если после очередного преобразования появляется нулевая форма, путь не может быть кратчайшим.

\medskip

Допустим теперь, что наш путь по двум сторонам макроплитки уровня $n$ можно привести к форме, содержащей некорректный участок $XYZ$. 

Пусть $P$ соединяет левый нижний угол с правым верхним. Тогда $T$ вкладывается в макроплитку уровня $n+1$ в качестве левой верхней, при этом $P$ является подпутем аналогичного пути из левого верхнего угла в правый верхний в более крупной макроплитке. Операцию вложения макроплитки в более крупную можно повторить любое число раз. Если $P$ соединяет правый нижний угол с левым верхним, то $T$ можно вложить в макроплитку уровня $n+1$ в качестве нижней подплитки. Пусть $P$ становится частью пути из левого нижнего угла в правый верхний. Далее вложения можно делать как для случая выше.

Итак $P$ является подпутем сколь угодно протяженного от него в обе стороны пути $P'$, идущего по сторонам достаточно большой макроплитки и соединяющего ее левый нижний и правый верхний углы. По лемме о некорректных участках, если $P$ приводится к пути, содержащему некорректный участок, то $P'$ приводится к нулю. Но это невозможно, как уже показано выше.

\end{proof}

{\bf Замечание.} Можно показать, что путь по двум соседним сторонам макроплитки после любого числа преобразований будет имеет форму $W_1UW_2$, где пути $W_1$, $W_2$ (которые могут быть просто одной вершиной) полностью лежат на границах макроплитки, а у подпути $U$ концы лежат на границе, а все остальные вершины лежат строго внутри.

\begin{lemma}[О расстоянии от края до выхода в подклейку] \label{pasting_distance}

1. Пусть вершина $X$ лежит на краю некоторой макроплитки $T$ принадлежащей комплексу $K$, вершина $Y$ принадлежит $T$, но не находится на ее границе, а из $Y$ существует выход в подклееную макроплитку уровня $n \geq 2$, не содержащую вершин на границе $T$. Тогда расстояние от $X$ до $Y$ в комплексе $K$ (длина кратчайшего пути по ребрам) не менее $2^{n-1}$.

\end{lemma}

\begin{proof}

Сначала заметим, что для любых вершин $X$ и $Y$ на комплексе, операция подклейки не меняет расстояния между ними понимаемом в смысле количества ребер графа в кратчайшем соединяющем пути. Действительно,
если $X$ и $Y$ после подклейки оказываются соединенными более коротким путем. Каждый его подпуть соединяющий две точки на границах подклееной макроплитки. может быть локально преобразован в подпуть не заходящий внутрь этой подклейки, согласно лемме о вынесении пути на границу. В этом случае более короткий путь существует и до подклейки.

\medskip

Вернемся к утверждению леммы. Проведем доказательство индукцией по уровню макроплитки $T$. Для уровней $2$ и $3$
подклеек внутри $T$ не проводится. Для уровня $4$ ближайшая вершина с подклееным ребром, в макроплитку, не содержащую $X$, лежит в середине любого ребра, выходящего на середину любого из граничных ребер. Расстояние от нее до границы $T$ будет равно $2$, при этом подклееная макроплитка будет иметь уровень $2$.

Пусть уровень $T$ равен $k>4$. 
Возьмем кратчайший путь $XWY$. Можно считать, что $W$ не содержит вершин на границе $T$.
Кроме того, можно считать, что $Y$ не лежит внутри какой-то из подплиток $T$, иначе $X$ можно взять на границе этой подплитки.
Значит, $Y$ обязательно лежит на одном из ребер, которые разделяют макроплитку $T$ на $6$ макроплиток при разбиении. Опять используя минимальность пути, замечаем, что $Y$ лежит в углу подклееной макроплитки. Из определения подклейки получаем, что уровень $Y$ был на $1$ меньше максимального в тот момент, когда эта подклейка производилась. Если после этого разбиений не было, то расстояние $XY$ не менее $2$, а также макроплитка при $Y$ -- второго уровня. Если после этого провели $l\geq 1$ разбиений, расстояние $XY$ увеличилось в $2^l$ раз, а уровень макроплитки при $Y$ возрос на $l$. Учитывая, что подклейки не меняют расстояний между точками, это завершает доказательство.

%\begin{corollary}
%Расстояние между двумя выходами в подклееные области не менее $2^{n-1}$, где $n$ уровень макроплитки, куда идет второй выход.
%Пусть путь $P$ имеет форму $XP_1YZ$, где  $P_1$ некоторый путь, выходы из вершин $X$ и $Y$ осуществляются по подклееным ребрам, а все остальные входы во все вершины, и выходы из них плоские. Уровень подклееной макроплитки, куда идет выход из $Y$ равен $n$. Тогда расстояние от $X$ до $Y$ не менее $2^{n-1}$.
%\end{corollary}

\end{proof}

%Таким образом, если в пути встречаются два ребра выходящие из вершин по подклееным ребрам, то они разделены, как минимум, одним плоским ребром.

\medskip

%Каждая подклееная макроплитка имеет ограниченные размеры, соответственно, путь, не имеющий нулевой формы и не возвращающийся из подклееной макроплитки будет иметь ограничения на длину.

\begin{lemma}[Об ограниченности пути, уходящего в подклееную часть] \label{pastingpath}

1. Пусть путь $P$ имеет вид $V_1AV_2BV_3CV_4$, где ребра из вершин $A$, $B$, $C$ ведут в подклееные области, а участки $V_1$, $V_2$, $V_3$, $V_4$ -- плоские. Тогда путь $P$ не может быть плоским (лежать в одной макроплитке).

2.Пусть путь $P$ начинается в вершине $X$, выход из $X$ идет в подклееную макроплитку уровня $n$. Кроме того, пусть $P$ не содержит ребер-выходов из подклееных плиток.
Тогда, если длина пути не менее  $2^{n+1}$, то он может быть приведен к нулевой форме.

%Пусть путь $P$ имеет форму $XP_1YP_2$, где $P_1Y$ -- плоский участок (лежащий в какой-то макроплитке) и выходящие из вершин $X$ и $Y$ ребра идет в подклееные плитки, и больше выходов в подклееную область или входов в вершину из подклееной области нет. Пусть также $P$ не может быть приведен к нулевой форме.
%Тогда длина $YP_2$ не более, чем удвоенная длина $XP_1Y$.
\end{lemma}

\begin{proof}

1. Допустим, путь проходит в подклееной макроплитке $T$. Тогда $T$ -- та макроплитка, куда уходит ребро из $C$. Тогда путь от $A$ до $C$ проходит по краю $T$. Рассмотрим макроплитку (или две макроплитки) к которой подклеивается $T$. Тогда путь от начала до $C$ проходит по этой одной или двум макроплиткам, при этом два раза проходит по ребру в подклееную область, что невозможно. (Один переход возможен, при переходе от одной макроплитке к другой.)

2. Пусть $X_i$ -- вершины $P$ где происходят входы в подклееные макроплитки $T_i$ и $X_1$ это $X$. Пусть их всего $s$. Пусть уровень $T_i$ равен $n_i$, $n_1=n$. Учитывая, что уровень подклееной макроплитки, как минимум, на $1$ меньше, чем уровень каждой из макроплиток, к которым она подклеивается, получаем $n_1>n_2>\dots >n_s$.

По лемме о выделении локального участка, длина $X_i X_{i+1}$ менее $2^{n_i}$. Тогда суммарная длина пути $P$ меньше чем $2^n+2^{n-1}+\dots 2^3+2^2<2^{n+1}$.

\end{proof}

%\section{Библиография}

%\addcontentsline{toc}{chapter}{\Numline {}Библиография}

%\markboth{}{Библиография}

%\begin{enumerate}

% Доклады, сделанные по материалам работы, приведены в конце списка литературы.

\end{document}